%% file: PAPER.arxiv.tex
\documentclass[12pt]{article} 
\usepackage{geometry}
\usepackage{amsmath}
\usepackage{amsfonts}
\usepackage{amssymb}
\usepackage[pdftex]{graphicx}
\DeclareGraphicsExtensions{.pdf,.png,.jpg}
\usepackage{algorithmic}
\usepackage{bm}
\usepackage{color} 
\usepackage{float}
\usepackage[color=white]{todonotes} 
\usepackage{soul} 

\usepackage{enumitem}

\usepackage[titletoc,title]{appendix}

\input{mathmacros.tex}
\begin{document}

\title{Passive array imaging  in random media}

\author{Liliana~Borcea and Ilker~Kocyigit
\thanks{L.~Borcea is with the Department of Mathematics, University of Michigan, Ann Arbor, 
MI 48109 USA,  e-mail: borcea@umich.edu.}
\thanks{I.~Kocyigit is with the Department of Mathematics, Dartmouth College, Hanover, NH 03755 USA, 
e-mail: ilker.kocyigit@dartmouth.edu}
}

\markboth{ }%
{Shell \MakeLowercase{\textit{Borcea and Kocyigit}}: Passive array imaging in  random media}
\maketitle

\begin{abstract}
We present a novel algorithm for high resolution coherent imaging of sound sources 
in random scattering media using time resolved measurements 
of the acoustic pressure at an array of receivers. 
The sound waves travel a long distance between the sources 
and receivers  so that they are significantly affected by scattering 
in the random medium.  We model the scattering effects by large random wavefront 
distortions, but the results extend to stronger effects, as long as the waves 
retain some coherence i.e., before the onset of wave diffusion. 
It is known that scattering in random  media can be mitigated in imaging using coherent interferometry (CINT). 
This method introduces a statistical stabilization in the image formation, at the cost of image blur. 
We show how to modify the CINT method in order to image 
wave sources that are too close to each other to be distinguished by CINT alone.
We introduce the algorithm from first principles and demonstrate its performance with numerical simulations.
\end{abstract}
\providecommand{\keywords}[1]{\textbf{\textit{Keywords---}} #1}

\section{Introduction}

\begin{figure}[!t]
\center
\includegraphics[width = 4.80in]{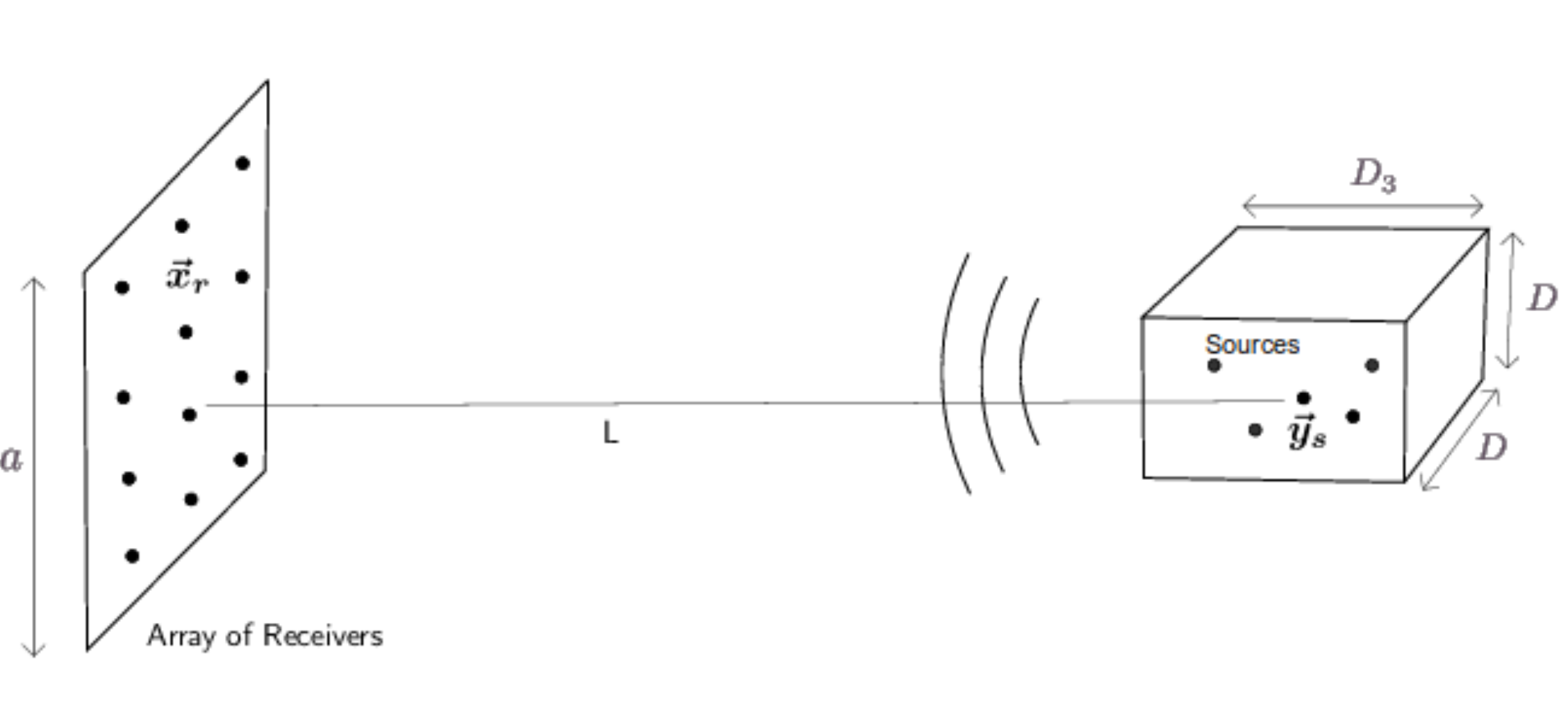}
\vspace{-0.15in}
\caption{Imaging setup with an array of receivers
  that is planar square of side $a$. The range direction is orthogonal
  to the array aperture. The unknown sources are  in the  imaging region, a rectangular prism with
  size $D_3$ in range and $D$ in cross-range. }
\label{fig:setup}
\end{figure}

Coherent array imaging is an important technology in radar \cite{radar}, sonar \cite{sonar}, 
seismic imaging \cite{seismic}, photoacoustic imaging \cite{PhotoAc1}, medical imaging with ultrasound \cite{medic}, and so on. 
We focus attention on passive array imaging, where a collection of $N_r$ receivers record 
waves generated by $N_s$ unknown sources. The receivers are located at  points $\vx_r \in \mathcal{A}$, for $r = 1, \ldots, N_r$,  
where $\mathcal{A}$ is the array aperture, assumed for convenience to be planar and square, of side $a$, as shown in Figure \ref{fig:setup}. 
The unknown sources are located at points $\vy_s \in \mathcal{D}$,  where $\mathcal{D}$ is the imaging region, a bounded set 
with center at distance $L$ from the array, in the direction orthogonal to the aperture, called the range direction. The coordinates 
in the plane orthogonal to this direction are called cross-range coordinates.  We let $\mathcal{D}$ be   a rectangular prism  
with square cross-section of side $D$ in the cross-range plane, satisfying $L \gg a > D$, and side $D_3$ in the range direction, 
satisfying $L \gg D_3$.  These scaling relations are typical in most imaging applications. 

The waves are modeled by the pressure $p(t,\vx)$, the solution of the equation 
\begin{equation}
\frac{1}{c^2(\vx)} \partial_t^2 p(t,\vx) - \Delta p(t,\vx) = \sum_{s = 1}^{N_s} f_s(t) \delta(\vx-\vy_s),
\label{eq:1}
\end{equation}
for $\vx \in \mathbb{R}^3$ and $t \in \mathbb{R}$,  with zero initial conditions on $p(t,\vx)$ and $\partial_t p(t,\vx)$
at time $t$ prior to the source excitation. The variable wave speed 
\begin{equation}
c(\vx) = c_o \left[ 1 + \sigma \mu\left(\frac{\vx}{\ell}\right)\right]^{-1/2}
\label{eq:2}
\end{equation}
models the heterogeneous medium consisting of a homogeneous background with constant wave speed $c_o$ and   
numerous weak inhomogeneities of size $O(\ell)$,   
commonly referred to as clutter. Because these inhomogeneities are unknown in imaging, they introduce uncertainty in the wave 
propagation, modeled in equation \eqref{eq:2} by the dimensionless random process $\mu$ of dimensionless argument.  We assume 
that $\mu$ is stationary,  and bounded almost surely. It has  zero mean,  and  autocorrelation 
\begin{equation}
{\cal{R}}_\mu(\vx-\vx') = \EE[ \mu(\vx)\mu(\vx')] = \exp \left[-\frac{|\vx-\vx'|^2}{2} \right],
\label{eq:Autocor}
\end{equation}
where $\EE[\cdot ]$ denotes expectation with respect to the distribution of $\mu$. The Gaussian 
expression of $\cal{R}_\mu$ is chosen for convenience, but the results extend to 
any integrable autocorrelation function. The scale $\ell$ is  called 
the correlation length and $\sigma \ll 1$ quantifies the small amplitude of the fluctuations of $c(\vx)$. 

The imaging problem is to determine the source locations $\{\vy_s\}_{1 \le s \le N_s}$ from the measurements $\{p(t,\vx_r)\}_{r=1,\ldots, N_r}$.

\subsection{Related work}  Each inhomogeneity in clutter is a  weak scatterer by itself since $\sigma \ll 1$, but  
cumulative scattering builds up over long ranges.  Mathematically, this manifests in the exponential decay of the coherent 
wave $\EE[p(t,\vx)]$ and the increase of the 
fluctuations $p(t,\vx) - \EE[p(t,\vx)]$. The range scale $\cal{S}$ of decay of the coherent wave is  the scattering mean free  path \cite{scattering}. 

When  $L < \cal{S}$,  the cumulative scattering effects  are negligible.  Much of the imaging literature considers this case, and  
coherent methods known as reverse time migration \cite{seismic}, matched filtering or backprojection 
\cite{radar,backproject} work well. They are based on the data model 
\begin{equation}
p(t,\vx_r) = p_o(t,\vx_r) + W(t,\vx_r), \quad r = 1, \ldots, N_r, 
\label{eq:3}
\end{equation}
where $p_o(t,\vx_r)$ is the solution of equation \eqref{eq:1}  with constant wave speed $c_o$, 
and $W(t,\vx_r)$ is additive noise with some statistics, assumed uncorrelated over the receivers. In the simplest form, the imaging function is 
\begin{equation}
\mathcal{J}(\vy) = \sum_{r = 1}^{N_r} p\left(\tau(\vx_r,\vy),\vx_r\right), \quad \vy \in \mathcal{D},
\label{eq:4}
\end{equation}
where 
\begin{equation}
\tau(\vx_r,\vy) = {|\vx_r-\vy|}/{c_o}
\label{eq:5}
\end{equation}
is the travel time from the imaging point $\vy$ to the receiver at $\vx_r$. 
The image  is robust to additive noise, and it  peaks in  the vicinity of the source locations, with cross-range 
resolution $O(\la_o L/a)$ and range resolution $O(c_o/B)$, where $\la_o = 2 \pi c_o/\om_o$ is the central  wavelength of the waves, calculated in terms of the central frequency $\om_o$,  and $B$ is the bandwidth 
of the source signals.  Better resolution can be achieved using  convex, sparsity promoting optimization, if the noise is not too strong, and 
the sources are separated by more than $\la_oL/a$ in cross-range and $c_o/B$ in range, see e.g.,
\cite{sparse, sparse1, sparse2}.

The imaging problem is much more difficult when $L \gg \cal{S}$, because the array measurements are significantly affected by scattering in clutter
and are no longer approximated by the model \eqref{eq:3}. If the range $L$ is so large that it exceeds the transport mean free path $\cal{T}$, 
which is the  scale that marks the onset of wave diffusion  \cite{scattering}, coherent imaging 
cannot succeed.  We assume an intermediate 
regime $\mathcal{T} > L \gg \mathcal{S}$, where coherent 
imaging is still possible. 
Such imaging must involve statistical stabilization with respect to the uncertainty of clutter, so that the estimates of the source locations 
are insensitive to the particular realization of the random medium (clutter) in which the imaging takes place.
 
Statistical stability  can be obtained with 
the coherent interferometric (CINT) approach \cite{CINT,CINT1}, which forms images using  the cross-correlations 
\begin{align}
{\cal C}(t,\tit,\vx_r,\vx_{r'}) =& \int_{-\infty}^\infty ds \, \Omega \, \Phi\big[\Omega(t-s)\big] \times \nonumber \\
&p\Big(s-\frac{\tit}{2},\vx_r\Big) p^\star\Big(s+\frac{\tit}{2},\vx_{r'}\Big), 
\label{eq:6}
\end{align}
where the star denotes complex conjugate. These are calculated around the time $t$, in a time window modeled by 
the function $\Phi$ of dimensionless argument and  $O(1)$ support. The  parameter $\Omega$  is adjustable and it 
should be similar to $\Omega_d$, the frequency offset over which the waves decorrelate in the random medium  \cite{scattering, CINT}.
Because the waves also decorrelate over directions of arrival, only cross-correlations at nearby receivers are useful. Thus, CINT
uses a spatial windowing function $\Psi$ of dimensionless argument and  $O(1)$ support to ensure that the receivers 
in \eqref{eq:6} are at distance $|\vx_r-\vx_{r'}| \le X$, with adjustable parameter that optimally equals  the 
decorrelation length \cite{CINT}. This scale is  proportional to the wavelength, but 
 typically $B \ll \om_o$, so the decorrelation length $X_d$ is approximately constant in the bandwidth.  The CINT imaging function is 
\begin{align}
{\cal J}(\vy) = &\sum_{r,r' = 1}^{N_r} \Psi \Big( \frac{|\vx_r-\vx_{r'}|}{X} \Big) \times \nonumber \\
&{\cal C}(\bar{\tau}(\vx_r,\vx_{r'},\vy),\widetilde{\tau}(\vx_r,\vx_{r'},\vy),\vx_r,\vx_{r'}),
\label{eq:7}
\end{align}
where 
\begin{align*}
\bar{\tau}(\vx_r,\vx_{r'},\vy) &= \frac{1}{2} \left[ \tau(\vx_r,\vy) + \tau(\vx_{r'},\vy)\right],
\\
\widetilde{\tau}(\vx_r,\vx_{r'},\vy) &= \tau(\vx_r,\vy) - \tau(\vx_{r'},\vy).
\end{align*}
It is statistically stable with respect to the realizations of the random medium if $\Omega \le \Omega_d \ll B$ and $X \le X_d \ll a$, 
meaning that its expectation near the peaks is much larger than the standard deviation \cite{CINT1}. Moreover,  these peaks 
are in the vicinity of the source locations with resolution $O(\la_o L/X)$ in cross-range and $O(c_o/\Omega)$ in 
range.  These resolution limits are similar to those in homogeneous media, except that the 
 aperture size $a$ and bandwidth $B$ are replaced  by the windowing parameters $X$ and 
$\Omega$. These are necessarily smaller than $a$ and $B$ to have statistical stability, so the images are  blurrier. 

It is shown in \cite{CINT_opt} that the resolution of CINT images may be improved using convex optimization. However, this requires 
detailed knowledge of the blurring kernel i.e., prior calibration. The deblurring  also works best when the sources are separated by distances 
larger than $\la_o L/X$ in cross-range and $c_o/\Omega$ in range, in the sense that if this is not so, there is no guarantee of unique 
recovery of the source locations. 
\subsection{Contributions}
In this paper we show how, by slightly modifying the CINT imaging function \eqref{eq:7}, it is possible 
to recover the unknown sources almost as well as in homogeneous media. Explicitly, we show that 
a collection of sources that are within a blurred peak of the CINT image and are separated by distances 
$O(\la_o L/a)$ in cross-range and $O(c_o/B)$ in range, can be estimated up to a translation in the support 
of the CINT peak. We introduce an algorithm that achieves this result,  motivate it from first principles and assess
its performance with numerical simulations. 

Note that although we restrict our study to passive array imaging, the results generalize easily to active arrays that probe the medium 
with waves and record the echoes,  in order to determine point-like scatterers with much larger reflectivity than $\sigma$, the reflectivity of
the clutter inhomogeneities. These scatterers are  secondary sources of waves, 
which emit signals proportional to the incident wave, so they can be viewed as the unknown sources considered here. This is obvious in the single scattering (Born) approximation, but it extends to multiple scattering as well, as  explained  in \cite{sparse3}.  Point-like scatterers 
play an important role in sonar and radar imaging because  corners of targets create stronger echoes
than other features. Thus,  the images are often 
 a constellation of peaks from which target features are to be extracted \cite{shaperadar}. Our algorithm can be used 
for this purpose in random media, because it recovers the relative location of the corner reflectors with the same resolution as 
in homogeneous media. 

The paper is organized as follows: In section \ref{sect:CINT} we introduce and analyze the CINT-like imaging function. 
In section \ref{sect:ALG} we use this function in the source reconstruction algorithm. We illustrate the performance of this algorithm with numerical simulations in section \ref{sect:NUM}, and conclude with a summary in section \ref{sect:SUM}.

\section{The CINT-like imaging function}
\label{sect:CINT}
We propose a simple modification of \eqref{eq:7}, where instead of searching for a single point $\vy \in {\cal D}$, we 
have two search points $\vy,\vy' \in {\cal D}$. The CINT-like imaging function 
\begin{align}
{\cal I}(\vy,\vy') = &\sum_{r,r' = 1}^{N_r} \Psi \Big( \frac{|\vx_r-\vx_{r'}|}{X} \Big) \times \nonumber \\
&\hspace{-0.5in}{\cal C}(\bar{\tau}(\vx_r,\vx_{r'},\vy,\vy'),\widetilde{\tau}(\vx_r,\vx_{r'},\vy,\vy'),\vx_r,\vx_{r'}),
\label{eq:10}
\end{align}
superposes the cross-correlations \eqref{eq:6} evaluated at the $\vy$, $\vy'$ dependent travel times
\begin{align}
\bar{\tau}(\vx_r,\vx_{r'},\vy,\vy') &= \frac{1}{2} \left[ \tau(\vx_r,\vy) + \tau(\vx_{r'},\vy')\right], \label{eq:defOTTau1}
\\
\widetilde{\tau}(\vx_r,\vx_{r'},\vy,\vy') &= \tau(\vx_r,\vy) - \tau(\vx_{r'},\vy').\label{eq:defOTTau2}
\end{align}
To explain why this is beneficial, we derive  below the expression of ${\cal I}(\vy,\vy')$ using a random travel time model 
\cite{RANDTT} that accounts for large, random wavefront distortions in random media,  as assumed in adaptive optics \cite{AO}. 
This model is convenient for the calculations and has been used in the analysis of CINT imaging in random media in \cite{CINT1,CINT_opt}.

As shown in the appendix \ref{sect:appA}, the calculations are based on the expression of the second statistical moments of the pressure waves, 
which are qualitatively similar to those in stronger scattering regimes \cite{RANDTT,CINTTH}. Thus, the results extend  verbatim to such regimes.

\subsection{Setup and the random travel time model}
\label{sect:assume}
We introduce here a few assumptions that simplify the calculations and lead to an explicit expression of  
\eqref{eq:10}. 

The first assumption is that the sources emit the same pulse 
\begin{equation}
f_s(t) = f(t) = \left(\frac{2}{\pi}\right)^{1/4} \sqrt{B} \exp[-i \om_o t - t B^2], 
\label{eq:13}
\end{equation}
modulated at central frequency $\om_o$, with Gaussian envelope normalized so that $\|f\|_2 = 1$. The Fourier transform 
of this pulse 
\begin{equation}
\hat f(\om) =  \left(\frac{\sqrt{2 \pi}}{B}\right)^{1/2} \exp \Big[ - \frac{(\om-\om_o)^2}{4 B^2} \Big],
\label{eq:14}
\end{equation}
is also a Gaussian, centered at the frequency $\om_o$ and with standard deviation proportional to $B$, called in an abuse of terminology
the bandwidth. We assume that $B \ll \om_o$. In general,   the sources will not 
emit the same signal and moreover, the signal may not be a pulse. 
The imaging algorithm in this paper applies to arbitrary  $f_s(t)$, that may even be noise-like, as long as for
nearby sources these signals are statistically correlated. This holds for example in active array imaging,  where  unknown 
scatterers act as secondary sources of waves and $f_s(t)$ are given by the convolution of the probing signal emitted by the array
and the Green's function that propagate the waves in the random medium to the scatterer locations.

The second assumption is that at any given frequency $\om$ in the support of \eqref{eq:14}, the wave propagation can be modeled 
by the Green's function 
\begin{equation}
\hat G(\om,\vx,\vy) = \frac{\exp \left\{ i \om \left[\tau(\vx,\vy) + \delta \tau(\vx,\vy) \right]\right\}}{4 
\pi |\vx-\vy|}
\label{eq:15}
\end{equation}
of Helmholtz's equation, where  
\begin{equation}
\delta \tau(\vx,\vy) = \frac{\sigma|\vx-\vy|}{2 c_o} \int_0^1 d u \, \mu \Big(\frac{(1-u)\vy + u \vx}{\ell} \Big).
\label{eq:16}
\end{equation}
This is the random travel time model and we refer to \cite{CINT1,RANDTT} for a detailed discussion of  its range of validity. Here it suffices to say 
that it holds when $\la_o \ll \ell \ll a \ll L$ and $\sigma$ is sufficiently small. 
We 
also recall from \cite[Lemma 3.1]{CINT1} that the random process \eqref{eq:16} is  approximately Gaussian, with mean zero and standard deviation $O(\sigma \sqrt{\ell L}/\la_o)$. 

The third assumption is that the size of the imaging region satisfies the scaling relations 
\begin{equation}
\frac{\la_o L}{X} \ll D \le a, \quad \frac{c_o}{\Omega} \ll D_3.
\label{eq:17p}
\end{equation}
The lower bounds in these relations are the CINT resolution limits, so this assumption 
ensures that the search domain is large enough to observe the focusing of the CINT 
imaging function \eqref{eq:7}. We also suppose that the aperture size $a$ and the 
the size of the imaging region are sufficiently small so that the rays connecting the sources and 
receivers are contained within a narrow cone of small opening angle and axis along the range direction. 
This assumption is described in more detail in the appendix \ref{sect:appA} and in technical terms it means that 
the waves are in a paraxial propagation regime.

Finally, to carry out explicit calculations, we take the Gaussian window functions
\begin{align}
\Phi(\Omega t) &= \exp\Big[-\frac{(\Omega t)^2}{2} \Big],
\label{eq:18} \\
\Psi \Big(\frac{|\vx_r-\vx_{r'}|}{X}\Big) &= 
\exp\Big[ - \frac{|\vx_r-\vx_{r'}|^2}{2 X^2} \Big],
\label{eq:19}
\end{align}
in definitions \eqref{eq:6} and \eqref{eq:10}.
We also suppose that the receivers are spaced at $O(\la_o)$ distances, so that $N_r = O( a^2/\la_o^2) \gg 1$. 
This allows us to 
approximate the sums over the receivers by integrals over the aperture $\cA$. To avoid specifying the aperture size in
these integrals, we use the Gaussian apodization  \[\exp \big[-|\vx|^2/(2 (a/6)^2)\big], \quad \vx \in \cA,\]
which is negligible outside the disk of radius $a/2$.

\subsection{The imaging kernel}
\label{sect:kernel}
For simplicity, we neglect additive noise in the calculations in this section, although  noise is considered in the numerical 
simulations  in section \ref{sect:NUM}.
As shown in the appendix \ref{sect:appA}, in the setup described in section \ref{sect:assume}, the imaging function \eqref{eq:10} can be written as 
\begin{equation}
{\cal I}(\vy,\vy') = \sum_{s,s'=1}^{N_s} {\cal K}(\vy,\vy',\vy_s,\vy_{s'}),
\label{eq:20}
\end{equation}
with kernel $\cal{K}$ that depends on the search points $\vy, \vy'$ and the locations $\vy_s$ and $\vy_{s'}$ of pairs of the unknown sources. 
To describe this kernel, consider the system of coordinates with origin at the center of the array and write $\vy = (\by,y_3)$, with two dimensional vector $\by$ in the cross-range plane and range coordinate $y_3$. Define also 
the center and difference cross-range vectors 
\begin{equation}
\cby = \frac{\by + \by'}{2}, \quad \tby = \by-\by',
\label{eq:21}
\end{equation}
and the center and difference range coordinates 
\begin{equation}
\cy_3 = \frac{ y_3 + y'_3}{2}, \quad \ty_3 = y_3-y_3'.
\label{eq:22}
\end{equation}
Similarly, we let $\vy_s = (\by_s,y_{s,3})$ and define 
\begin{align}
\cby_{ss'} &= \frac{\by_s + \by_{s'}}{2}, \quad \tby_{ss'} = \by_s-\by_{s'},
\label{eq:21p} \\
\cy_{ss',3} &= \frac{ y_{s,3} + y_{s',3}}{2}, \quad \ty_{ss',3} = y_{s,3}-y_{s',3}, 
\label{eq:22p}
\end{align}
for all $s, s' = 1, \ldots, N_s.$

The kernel satisfies 
\begin{align}
|{\cal K}(\vy,\vy',\vy_s,&\vy_{s'})| \sim \exp \left\{ - \frac{|\tby_{ss'}|^2}{2 \gamma X_d^2}  - \frac{|\ty_{ss',3}-\ty_3|^2}{2 (c_o/B)^2}  
\right. \nonumber \\
& 
          -\frac{|\tby_{ss'}-\tby|^2}{2 [\gamma_1 L/ (k_o a)]^2} 
          - \frac{|\cby_{ss'}-\cby|^2}{2 [L/(k_o X_e)]^2} \nonumber \\
&\left. -\frac{(|(\cby_{ss'},\cy_{ss',3})|-|(\cby,\cy_3)|)^2}{2 (c_o/\Omega_e)^2} \right\}, \label{eq:kernel}
\end{align}
where  $\sim$ means of the order of, up to a multiplicative constant. We refer to the appendix \ref{sect:appA} for the detailed expression 
of ${\cal K}$, not just its absolute value. 
In \eqref{eq:kernel} we introduced the wavenumber 
\[k_o = \om_o/c_o = 2 \pi/\la_o\] and the positive coefficients $\gamma, \gamma_1 \ge  O(1)$. We also denote by  $\Omega_e$ and $X_e$ the frequency and length scales defined by 
\begin{align}
\frac{1}{\Omega_e^2} &= \frac{1}{\Omega^2} + \frac{1}{\Omega_d^2} + \frac{1}{4 B^2},
\label{eq:23} \\
\frac{1}{X_e^2} &= \frac{1}{X^2} + \frac{1}{X_d^2} + \frac{1}{4 (a/6)^2}.
\label{eq:24}
\end{align}
Since $B \gg \Omega_d$ and $a \gg X_d$, the optimal windowing choice $\Omega = O(\Omega_d)$ and $X = O(X_d)$  gives 
\begin{equation}
\Omega_e = O(\Omega_d) \ll B, \quad X_e = O(X_d) \ll  a. 
\label{eq:25}
\end{equation}

The expression \eqref{eq:kernel}  says that the imaging function peaks when the center $\frac{1}{2}(\vy + \vy)$
of the imaging points  is in the vicinity of $\frac{1}{2}(\vy_s + \vy_{s'})$, for some pair $s,s'$ of source indexes. The radius of this vicinity 
is    $O(\la_o L/X_e)$ in the cross-range plane and  $O(c_o/\Omega_e)$ in the range direction. This is  the same as the focusing 
of the CINT image \eqref{eq:7}. The new observation is that we can get much better estimates of the source offsets $\vy_s-\vy_{s'}$, with the same resolution 
as in the homogeneous medium i.e., $O(\la_o L/a)$ in cross-range and $O(c_o/B)$ in range. 
\section{The imaging algorithm}
\label{sect:ALG}
The imaging algorithms consists of the following steps:

\vspace{0.15in}
\textbf{Step 1:} Calculate the CINT image \eqref{eq:7}, which is the same as ${\cI}(\vy,\vy)$,  and identify its peaks in the search domain ${\cal D}$. 
The unknown sources lie in the support of these peaks, but they cannot be identified due to the poor resolution: $O(\la_o L/X_e)$ in cross-range and $O(c_o/\Omega_e)$ in range.  To reduce the computations, it suffices to form the CINT image on a coarse mesh with pixel size similar to  these resolution limits. 

\vspace{0.15in}
\textbf{Step 2:} Let $\vz_o$ be the location of the center of a CINT peak. Suppose that there are $n_{s} \le N_s$ sources within this peak and 
denote by ${\cal Y}$ the set of their locations. We can only expect to determine these 
locations up to an overall translation, so we set $\vz_o$ as one point in the constellation of $n_{s}$ sources.
To estimate the other locations, relative to $\vz_o$, calculate ${\cI}(\vz_0,\vy)$ for $\vy$ in the support of the CINT peak, on a refined imaging mesh with pixel size $O(\la_o L/a)$ in cross-range and $O(c_o/B)$ in range. These are the resolution limits for the source  offsets  in equation \eqref{eq:kernel}.  

\vspace{0.15in}
\textbf{Step 3:} 
Identify the peaks of ${\cI}(\vz_0,\vy)$, which are the points $\vz_j$ satisfying 
\begin{equation}
\vz_j-\vz_o \in \mathcal{E}({\cal Y}), \label{eq:26}
\end{equation}
where 
\begin{equation}
\mathcal{E}({\cal Y}) = \left\{ \vy_s-\vy_{s'}: ~ \vy_s, \vy_{s'} \in {\cal Y}, ~ \vy_s \ne \vy_{s'} \right\}.
\label{eq:26p}
\end{equation}
The set $\mathcal{E}({\cal Y})$ has  cardinality $n_s(n_s-1)$  in most cases,  
where the offset vectors  $\vy_s-\vy_{s'}$ are distinct for different pairs $(s,s')$, with $s, s' = 1, \ldots, n_s$ and $s \ne s'$. Thus, 
${\cI}(\vz_0,\vy)$ is expected to have $N_z = n_s(n_s-1)$ peaks. This count reflects that if 
$\vec{{\bf e}} \in \mathcal{E}({\cal Y})$, then $-\vec{{\bf e}} \in \mathcal{E}({\cal Y})$, as well. 
However,  there are special, unlikely cases, where different source pairs give the same offset vectors. Thus,  in general,
\[
N_z \le n_s(n_s-1).\]

\vspace{0.15in}
\textbf{Step 4:} From the $N_z$ peaks of ${\cI}(\vz_0,\vy)$ estimate the set \eqref{eq:26p} by 
\begin{equation}
{\cal E}^{\rm{est}} = \{ \vz_j - \vz_o, \quad j = 1, \ldots, N_z \}.
\label{eq:27}
\end{equation}
Use this set to determine the constellation of sources. For our purpose, it suffices to use the exhaustive search algorithm given below, which is  not optimal 
in terms of computational cost.  The output of this algorithm is a 
set ${\cal Y}^{\rm{est}}$ of points $\vy^{\rm{est}}_s$,  so that 
\begin{equation}
\mathcal{E}({\cal Y}^{\rm{est}}) = {\cal E}^{\rm{est}}.
\label{eq:28}
\end{equation}
Here $\mathcal{E}({\cal Y}^{\rm{est}})$ is defined as in \eqref{eq:26p}, with ${\cal Y}$ replaced by ${\cal Y}^{\rm est}$. 
The vectors in ${\cal Y}^{est}$ are the estimates of the source locations, up to the translation defined by fixing one source at $\vz_o$ and the reflection about $\vz_o$. 

\vspace{0.15in}
\begin{minipage}{\textwidth}

\begin{algorithmic}
\STATE{\textbf{Algorithm}$(E,Y)$ }
\STATE{\emph{Input:} The sets $E$, and $Y$.}
\STATE{\emph{Output:} Empty set or a non-empty set $\mathcal{Y}^{est}$.}
\vspace{0.05in}
\STATE{  1.~~\textbf{If}~~$\mathcal{E}\left(Y\right)=\mathcal{E}^{est}$~\textbf{then}~ } 
\STATE{		2.~~~~~~return~$\mathcal{Y}^{est}=Y$~~ }
\STATE{		3.~~\textbf{End-If}~ }
\vspace{0.05in}
\STATE{		4.~~\textbf{While}~$E\not=\emptyset$~~\textbf{do} } 
\STATE{		5.~~~~~~Select~the~first~vector~$\vec{\bf{e}}$~and~set~$E=E\backslash\left\{ \vec{\bf{e}}\right\} $~~ }
\STATE{		6.~~~~~~\textbf{Let}~$\vy=\vz_{0}+\vec{\bf{e}}$~  } 
\STATE{		7.~~~~~~\textbf{If}~$\left\{ \pm\left(Y-\vy\right)\right\} \subset\mathcal{E}^{est}$~~\textbf{then}~ } 
\STATE{		8.~~~~~~~~~~${Y}'=\textbf{Algorithm}(E,Y\cup\left\{ \vy\right\} )$~ } 
\STATE{		9.~~~~~~~~~~\textbf{If}~${Y}'\not=\emptyset$~\textbf{then}~  } 
\STATE{		10.~~~~~~~~~~~~~return~$\mathcal{Y}^{est}={Y}'$~~ } 
\STATE{		11.~~~~~~~~~\textbf{End-If}  }	
\STATE{		12.~~~~~\textbf{End-If}  }
\STATE{		13.~~\textbf{End-While}   }
\vspace{0.05in}
\STATE{     14.~~return~$\emptyset$} 
\end{algorithmic}

\end{minipage}

\subsection{Discussion}
The first call of this recursive algorithm is made with the inputs $E = \mathcal{E}^{est}$ and $Y = \{\vz_o\}$. 
When the algorithm outputs the empty set $\emptyset$, the search has failed.  We show in  the appendix 
\ref{sect:appB} that, when noise is not an issue i.e., the offset vectors in the set  ${\cal E}^{\rm{est}}$ are the same as 
those in $ \mathcal{E}({\cal Y})$,  
the output  of the algorithm is necessarily a non empty set  $\mathcal{Y}^{est}$ satisfying \eqref{eq:28}.
In practice, the testing of the equalities and the inclusion at lines 1 and 7  can be done up to some tolerance. In our numerical  simulations we consider two vectors to be same if  their difference has cross-range and range components that are smaller than the pixel size, in absolute value.

Note that the expression \eqref{eq:kernel} of the imaging kernel indicates that only sources at 
cross-range offsets $|\by_s-\by_s'| \le O(\sqrt{\gamma} X_d)$ contribute to the image ${\cI}(\vz_0,\vy)$ 
calculated at step 2. This distance is at least $O(\la_o L/X_e)$ in our case, so all the sources supported in 
the CINT peak should contribute to ${\cI}(\vz_0,\vy)$.  In other scattering regimes it may be that $X_d \ll \la_o L/X_e$, 
so the support of the CINT peak may be divided in smaller sets at step 2. The remainder of the algorithm above 
can be used separately for each such set.  

We already stated that we can only hope to determine the set ${\cal Y}$ of source locations up to an overall 
translation, fixed by the starting point $\vz_o$ and up to a reflection. This is due to the fact that 
the set $\mathcal{Y}^{\rm{ref}}$, defined as the reflection of $\mathcal{Y}$ with respect to a fixed point, satisfies 
$
\mathcal{E}({\cal Y})  = \mathcal{E}({\cal Y}^{\rm{ref}})  .
$

\subsection{Relation to localization} 
The estimation of  the source locations from the set $\mathcal{E}^{est}$ of  offset vectors  is somewhat related to 
the network localization problem 
\cite{networkLocal1}, \cite{networkLocal2}:
Let ${\cal N} = \left\{ \bn_1,\bn_2,\dots,\bn_m \right\}  $ denote the unknown set of nodes of a network. Determine ${\cal N}$ from knowledge of a  non-empty subset $\cal{ B} \subset {\cal N}$ of so-called  beacon nodes and 
the distance map $\delta: {\cal N} \times {\cal N} \to \mathbb{R}^+$, defined by  $\delta (i,j) = \|\bn_i-\bn_j\|$, for $i,j  = 1, \ldots, m$. 

Our problem is different as follows: (1) The data are the offset vectors $\vy_j-\vy_{j'}$ and not
just their norm. (2) We do not have access to the distance map $\delta$, we only know its image. 
That is to say, for any offset vector $\vec{\bf e} \in {\cal E}^{\rm{est}}$, we do not know the pair $(i,j)$ of sources 
that give $\vy_i-\vy_{j} = \vec{\bf e}$. (3) As explained at Step 3 in section \ref{sect:ALG}, in general, we do not know the number of  sources. Only under the 
additional assumption that there is a unique pair of sources that gives an offset vector $\vec{\bf e} \in {\cal E}^{\rm{est}}$,
we can determine the number of sources from the cardinality of the set ${\cal E}^{\rm{est}}$.
 
\section{Numerical results}
\label{sect:NUM}

\begin{figure}[t]
\center\includegraphics[width=4.4in,height=0.7in]{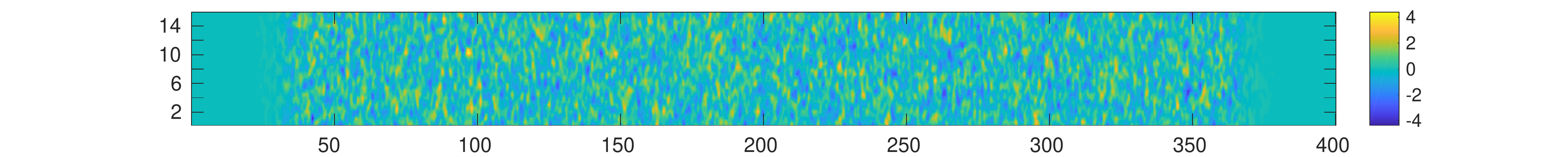}
\caption{Display of the realization of the fluctuations $\mu$ used in the numerical results.  The abscissa 
is range in units of $\ell$ and the ordinate is cross-range in units of $\ell$.}
\label{fig:2}
\end{figure}

To minimize the computational cost, we present imaging results in two dimensions, in the plane define by the 
range axis and one cross-range direction. The array cross-section in this plane is the line segment $[-a/2,a/2]$.

\subsection{Simulations setup}
The data are obtained using the model 
\begin{equation}
p(t,\vx_r) = \int_{-\infty}^\infty \frac{d \om}{2 \pi} e^{-i \om t} \left[\hat p(\om,\vx_r) + \hat W(\om,\vx_r)\right],
\label{eq:N1}
\end{equation}
where $\hat W$ denotes additive noise and 
\begin{equation}
\hat p(\om,\vx_r) = \hat f(\om) \sum_{s=1}^{N_s} \hat G(\om,\vx_r,\vy_s)
\label{eq:N2}
\end{equation}
is the solution of Helmholtz's equation in the cluttered medium. The noise $\hat W$ is complex Gaussian, uncorrelated 
over the receivers and frequencies, with mean zero and standard deviation $5\%$ of the maximum absolute value of \eqref{eq:N2}. The Green's function
$\hat G$  is calculated using definitions
\eqref{eq:15}-\eqref{eq:16}, in the realization of the random medium displayed in Figure \ref{fig:2}. This realization is generated 
using random Fourier series \cite{Fourier}, for the autocorrelation \eqref{eq:Autocor}.

All the length scales are relative to the correlation length $\ell$. The central wavelength is $\la_o = 1.1 \cdot 10^{-5} \ell$, the array aperture size is 
$a = 16 \ell$ and the range scale is $L = 800\ell$. The frequencies are scaled with respect to $\om_o$ and the bandwidth is $B = \om_o/5$.
The strength of the fluctuations  is $\sigma = 2 \cdot 10^{-6}$ and the decorrelation frequency and length defined in \eqref{eq:A4}
are  $\Omega_d = 0.039 \om_o$ and $X_d = 0.068 \ell$. The window parameters are $X = X_d/3$ and $\Omega = \Omega_d/3$, so 
definitions \eqref{eq:23}--\eqref{eq:24} give $X_e = 0.0214 \ell$ and $\Omega_e = 0.0124 \om_o$ and 
the decay scales in the expression of the kernel are 
\begin{align*}
\frac{L}{k_o X_e} = 0.0654 \ell, ~~ \frac{c_o}{\Omega_e} = 1.4 \cdot 10^{-4} \ell, \\
\frac{L}{k_o a} = 8.8 \cdot 10^{-5} \ell, ~ ~ \frac{c_o}{B} = 8.8 \cdot 10^{-6} \ell. 
\end{align*}
The coefficient $\gamma$ in \eqref{eq:kernel} is 
\[
\gamma = \frac{4 X_d^2}{4 X_d^2 - X_e^2} = 1.025
\]
and $\gamma_1  \ge 6 \sqrt{2}$.  

The images displayed in the next section are calculated on an imaging mesh with 
pixel size $3.92 L/(k_o a)$ in cross-range  and 
$0.537 c_o/B$ in range.  The axes in the plots are in units of the CINT resolution limits.  

\subsection{Source reconstructions}

\begin{figure}[H]
	\centering \includegraphics[width=0.31\textwidth]{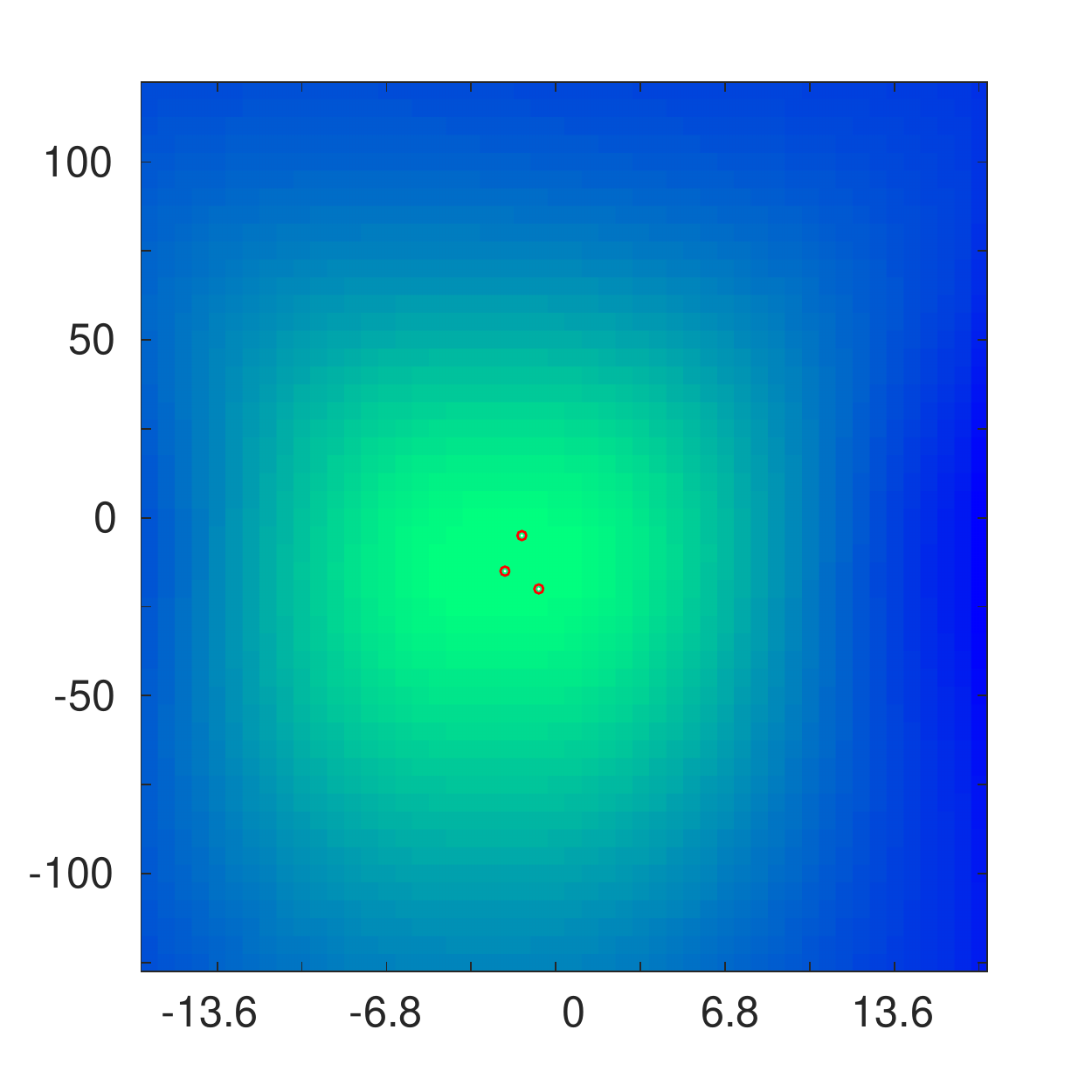}
	\hspace{-0.1in}\includegraphics[width=0.31\textwidth]{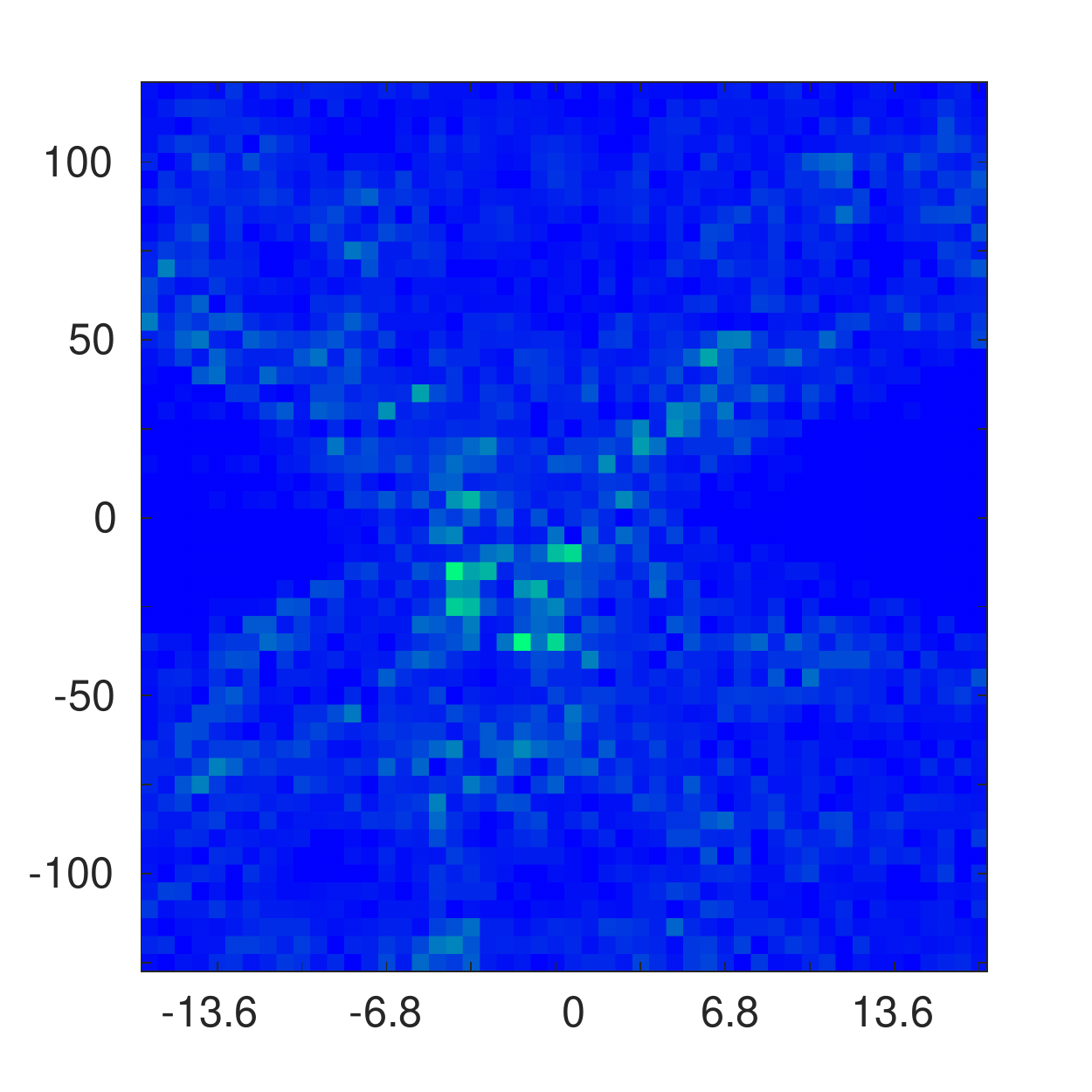} \\
	\includegraphics[width=0.31\textwidth]{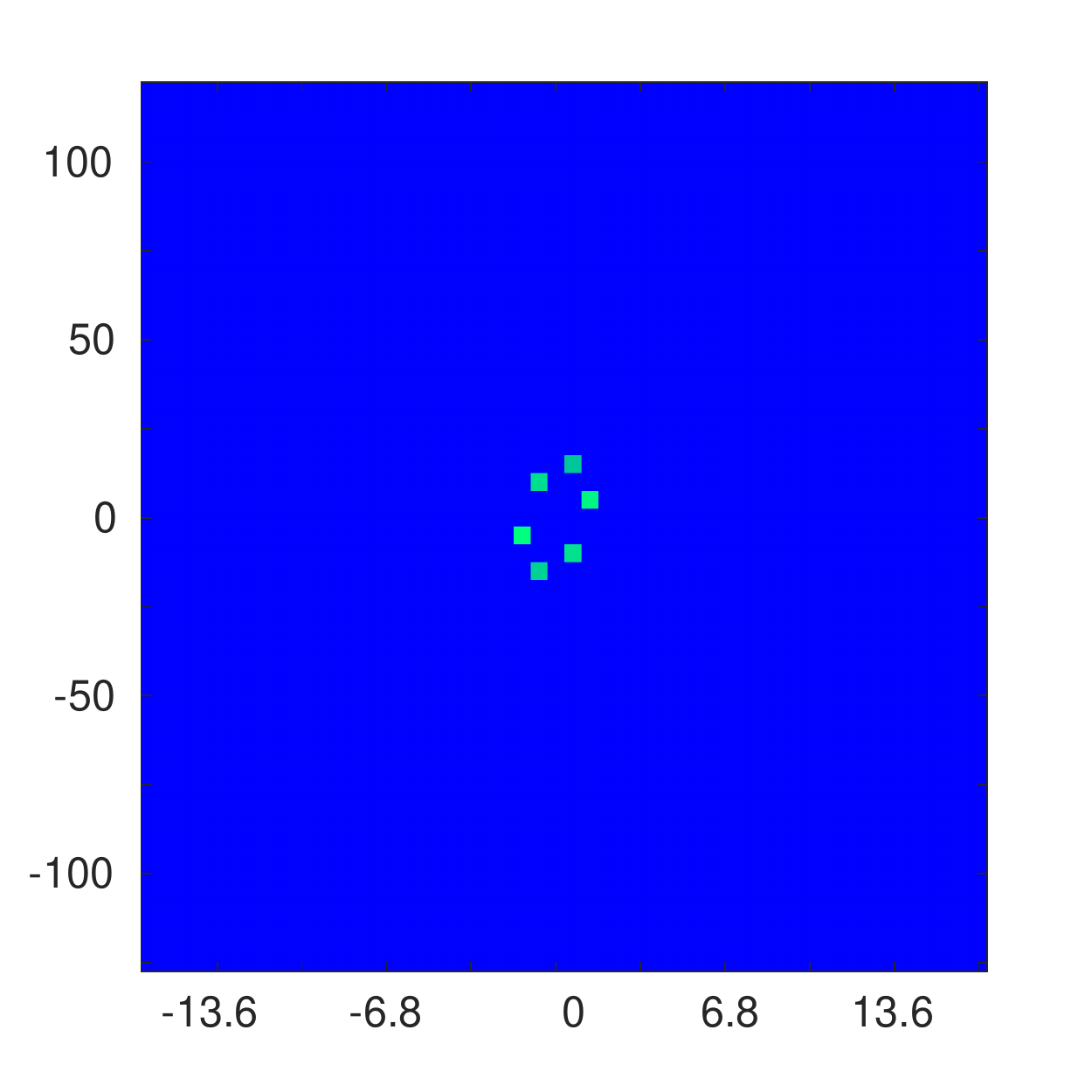}
	\hspace{-0.1in}\includegraphics[width=0.31\textwidth]{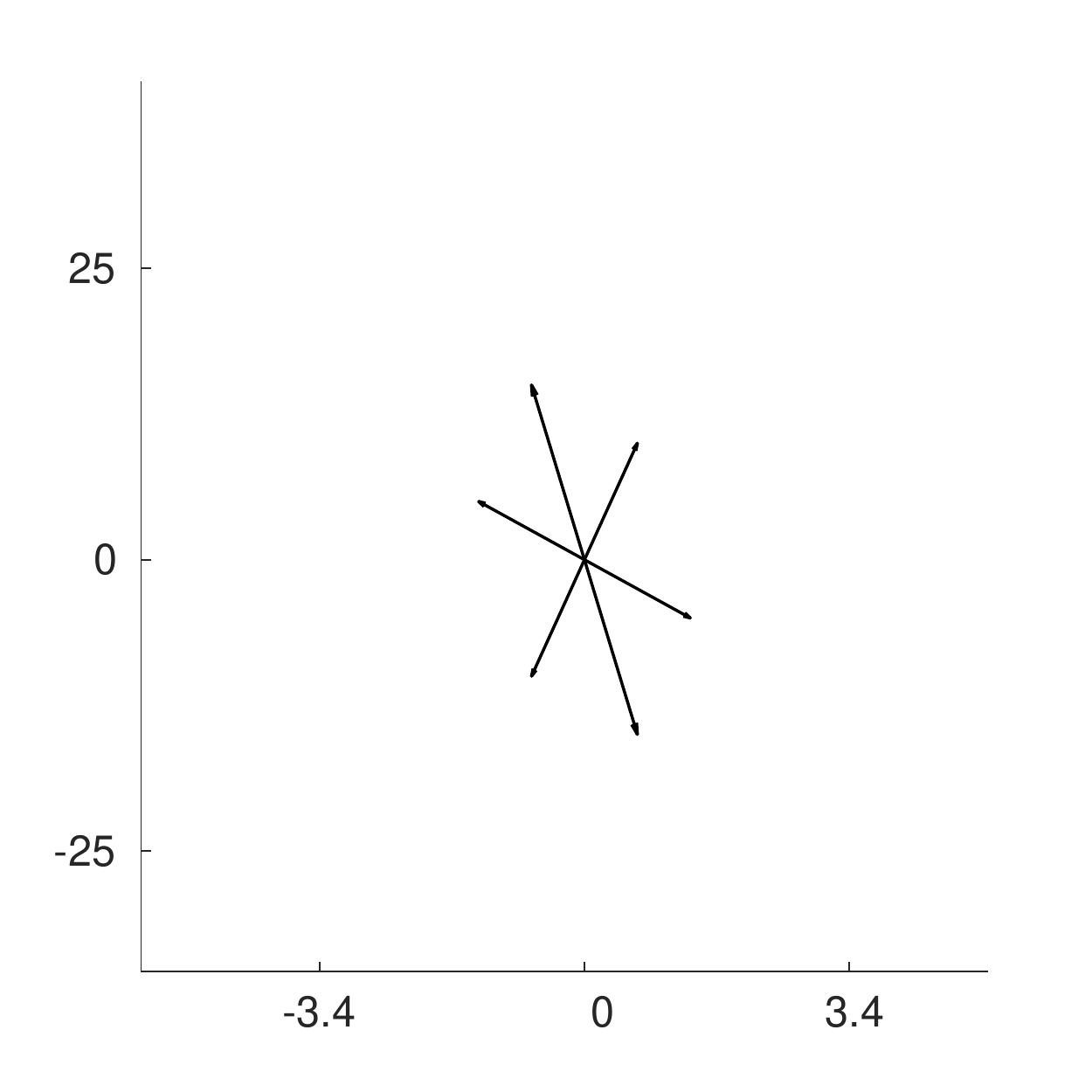} \\
	\includegraphics[width=0.31\textwidth]{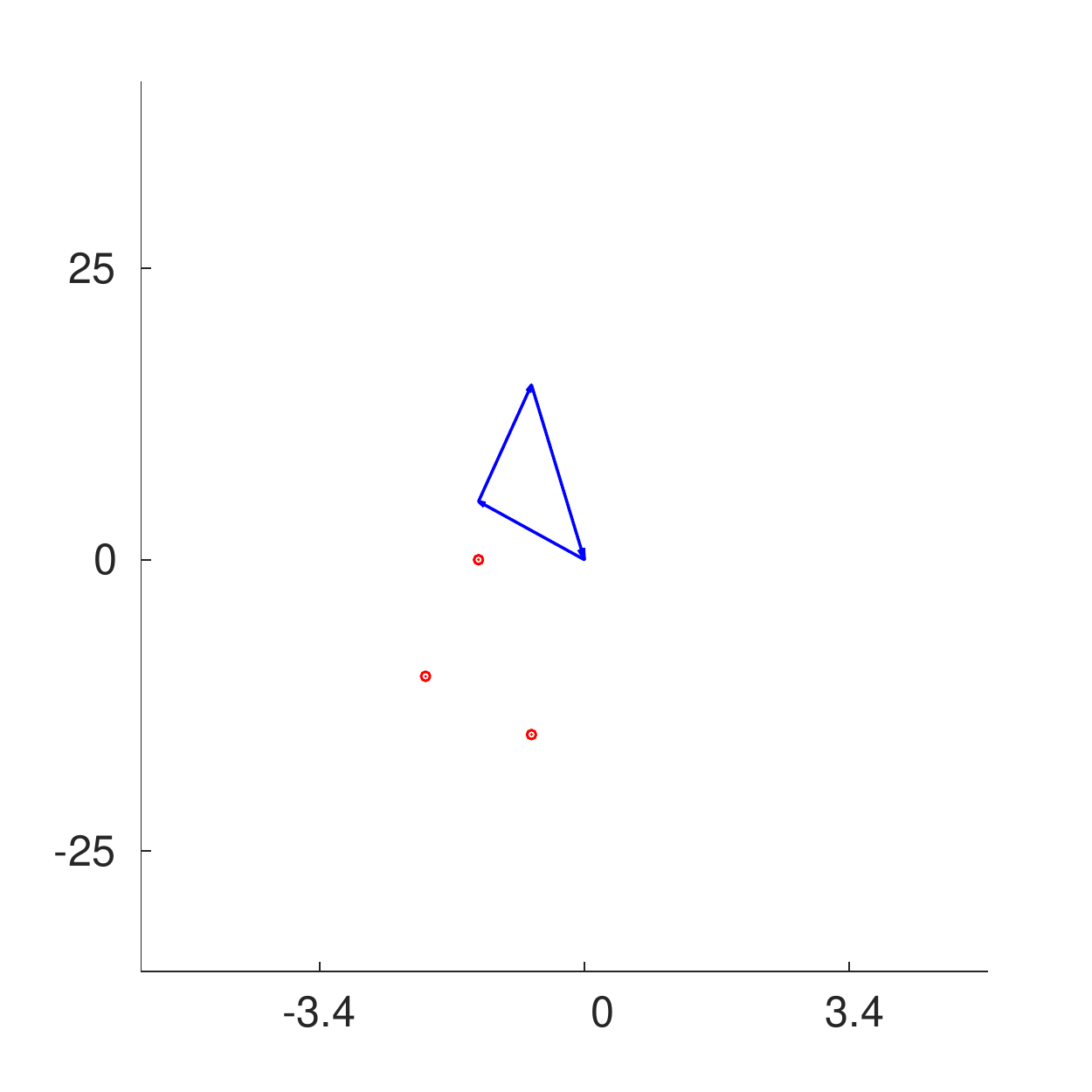}
	\hspace{-0.1in}\includegraphics[width=0.31\textwidth]{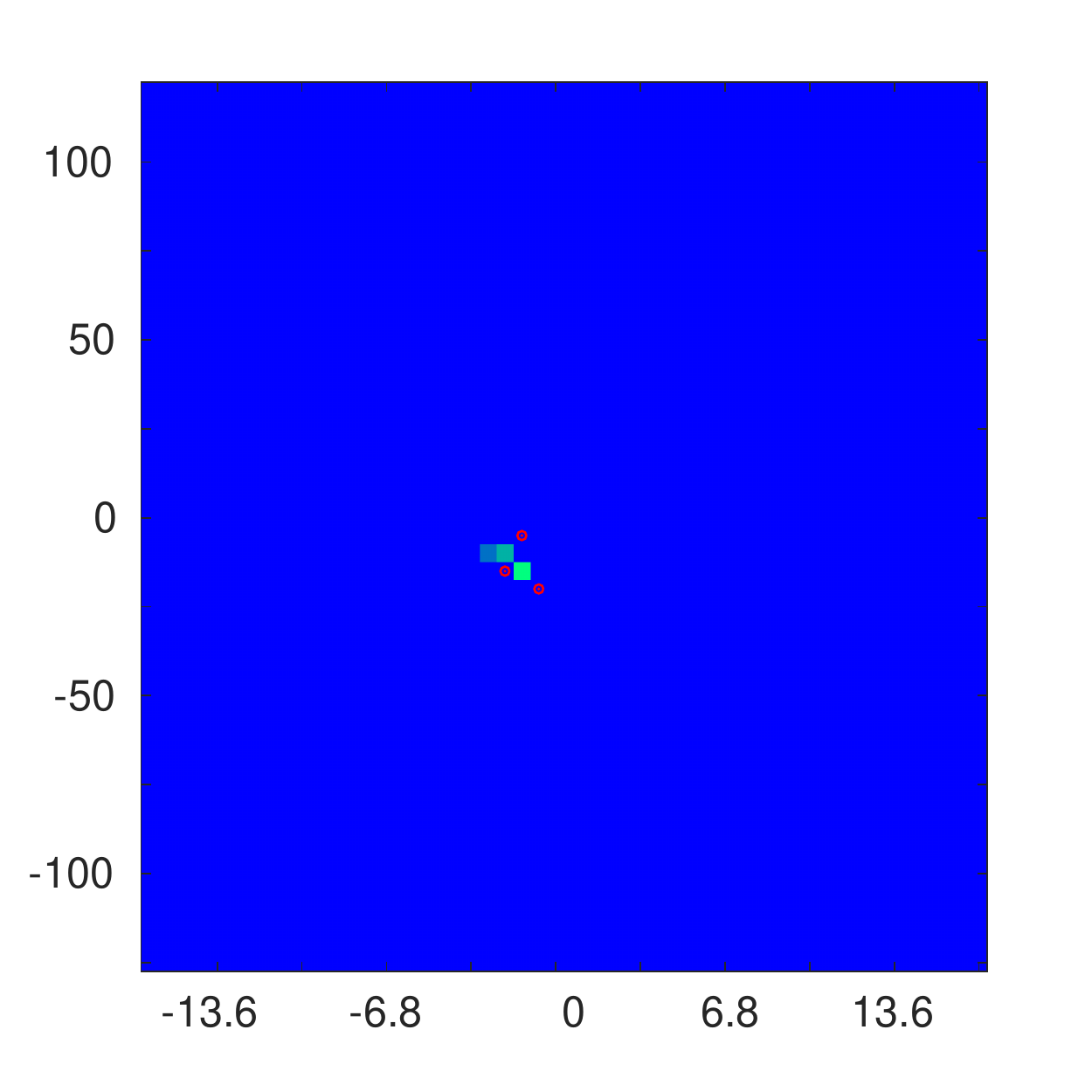}
	\caption{Top row: CINT image ${\cal I}(\vy,\vy)$ of three sources shown in red (left plot) and the conventional image
		\eqref{eq:4} (right plot).  Middle row:
		The image ${\cal I}(\vz_0,\vy)$ (left plot) and the  offset vectors that define the set \eqref{eq:27} (right plot).
		Bottom row:   Reconstruction of the three sources (left plot)  and reconstruction by deblurring the CINT image using convex optimization (right plot). The abscissa is range offset with respect to the center location $\vz_o$ of the CINT peak,     in units of $L/(k_o a)$.
		The ordinate is cross-range offset with respect to $\vz_o$, in units of $c_o/B$. 
	}
	\label{fig:3}
\end{figure}

We present in Figure \ref{fig:3} the results for three nearby sources. As expected, the CINT image ${\cI}(\vy,\vy)$ displayed 
in the top left plot has a blurry peak, and cannot distinguish the sources. The conventional image displayed in the top right plot  is calculated using definition \eqref{eq:4}. It  has many spurious peaks and  our simulations show that these peaks change unpredictably from one realization of the random 
process $\mu$ to another. This statistical instability is expected, because the data are incoherent in our regime. 

In the left plot of the middle row we display the image ${\cI}(\vz_o,\vy)$, with $\vz_o$ at the center of the CINT peak. 
Because we have $n_s = 3$ sources supported in this peak, we observe $N_z = 6$ peaks $\vz_j$, for $j = 1, \ldots, 6$. 
These define the set ${\cal E}^{\rm{est}}$ defined in \eqref{eq:27}, with offset vectors displayed in the right plot of the middle 
row. Note that for each vector $\vec{\bf e}$ in ${\cal E}^{\rm{est}}$ we also have the vector $-\vec{\bf e}.$ The reconstruction of the 
sources using the algorithm described in section \ref{sect:ALG} is shown in the bottom left plot.  The reconstruction is exact
up to the translation by the vector $\vz_o$, where by exact we mean with error that is smaller than the pixel size.
For comparison, we also show in the bottom right plot the reconstruction obtained with the debluring algorithm 
introduced in \cite{CINT_opt}. This algorithm is guaranteed to give a good reconstruction of the sources when they are further apart than  the CINT resolution limits. In this simulation the sources are much closer to each other so the results are worse than those in the bottom left plot.  As 
predicted by the theory in \cite{CINT_opt}, the reconstruction in the bottom right plot is peaked near the source locations, but it does 
not show three distinct sources.

\begin{figure}[h]
	\centering \includegraphics[width=0.31\textwidth]{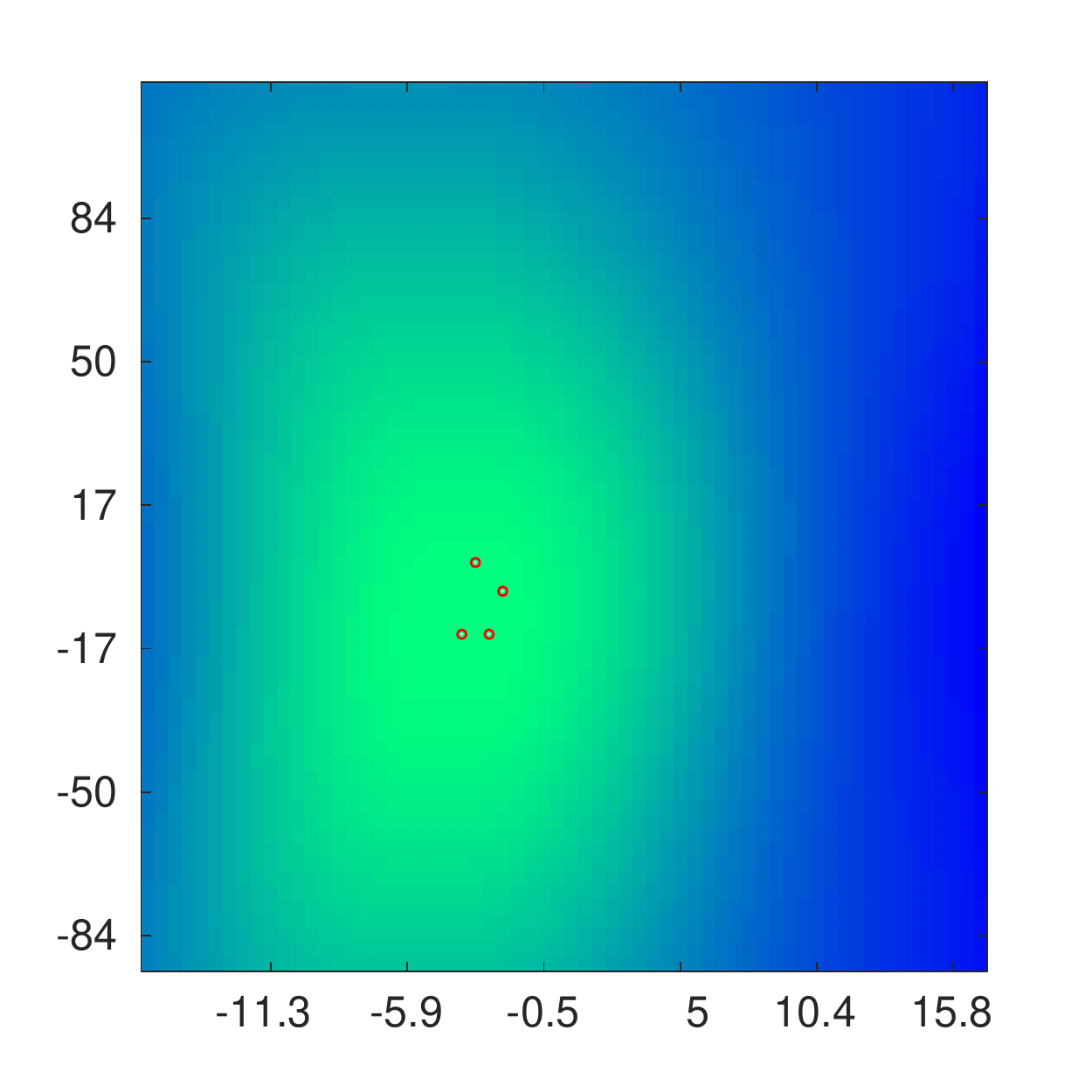}
	\hspace{-0.1in}\includegraphics[width=0.31\textwidth]{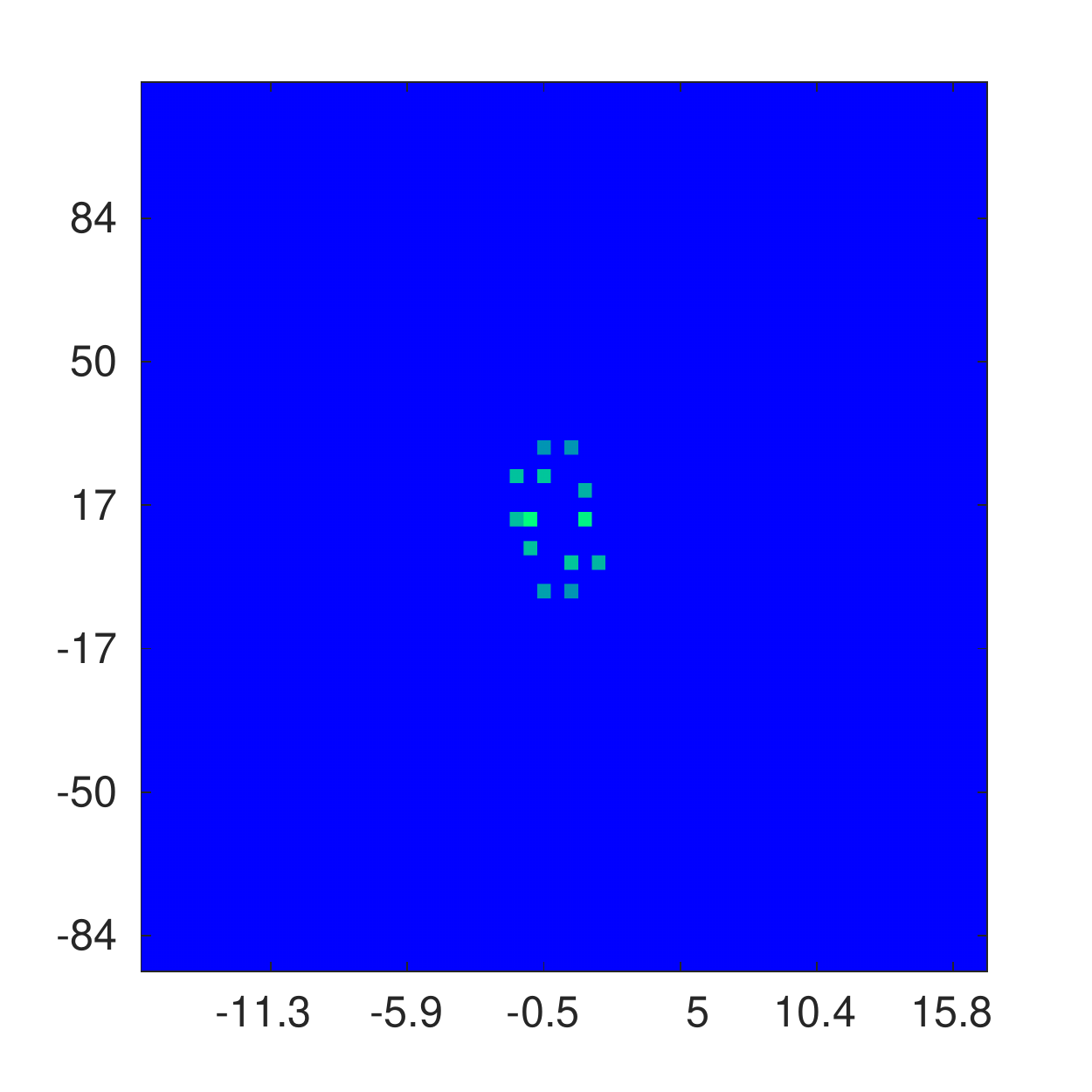}\\
	\centering \includegraphics[width=0.31\textwidth]{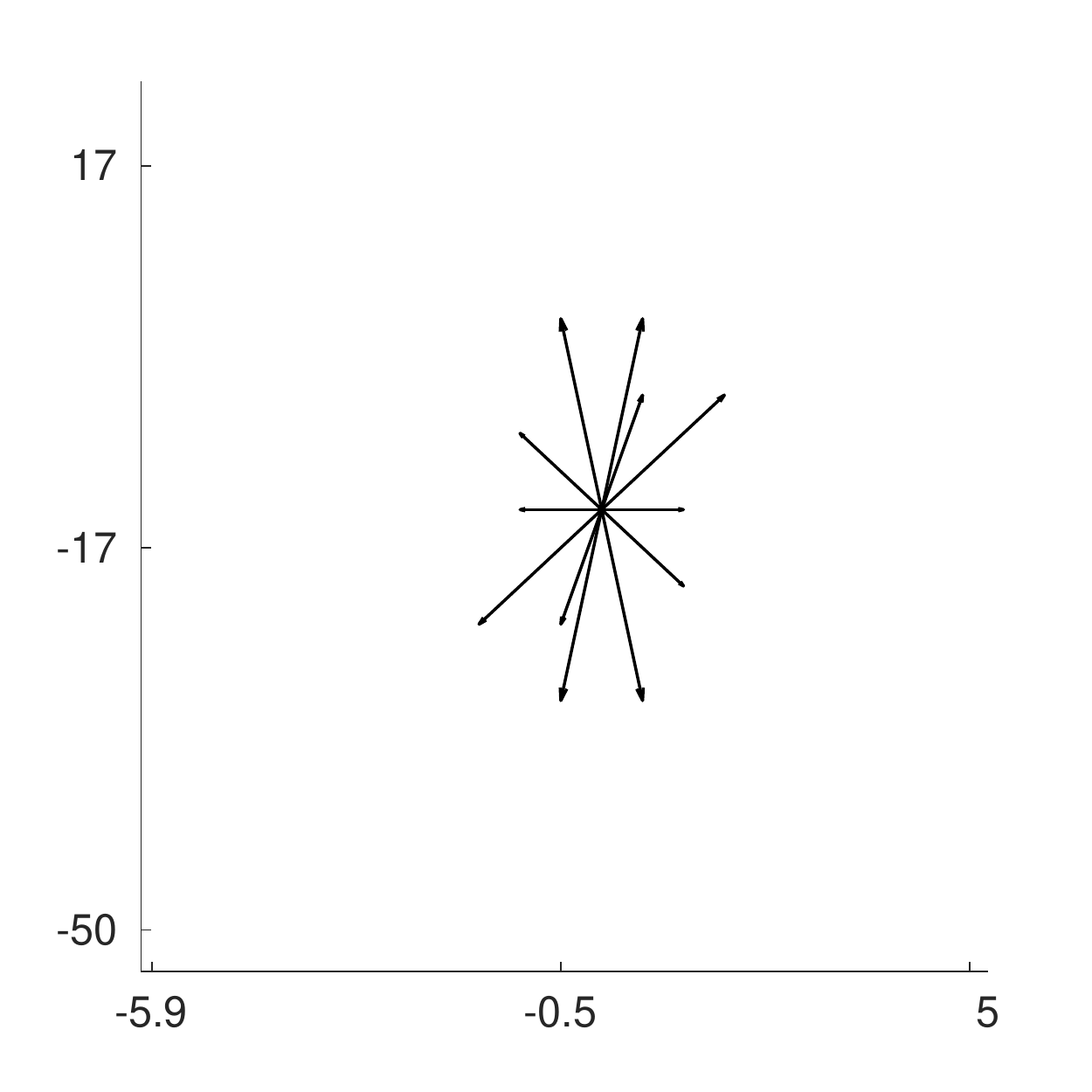}
	\hspace{-0.1in}\includegraphics[width=0.31\textwidth]{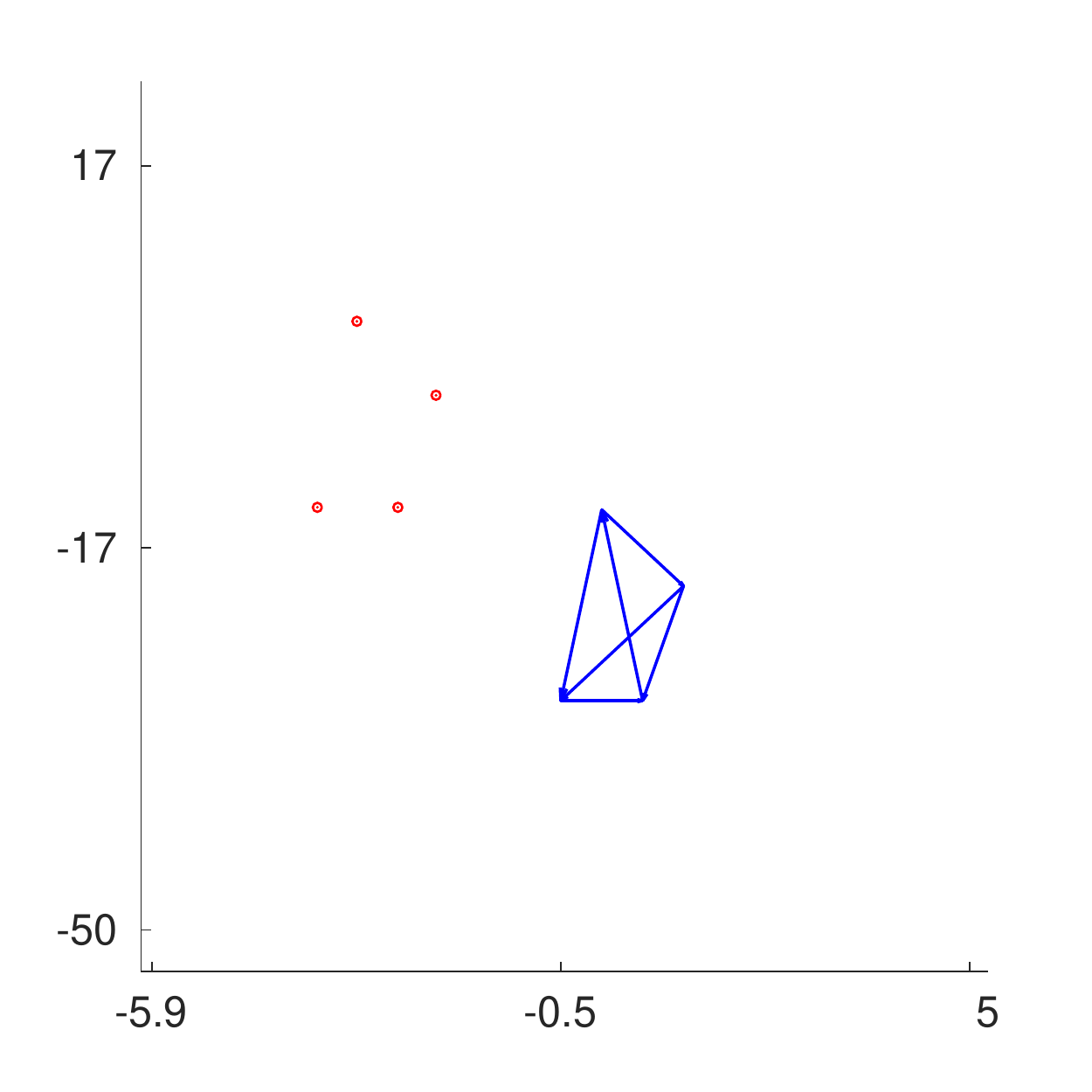}
	\caption{Top row: Left:  CINT image ${\cal I}(\vy,\vy)$ of four sources shown in red. Right: The image ${\cal I}(\vz_0,\vy)$.
		Bottom row:  Left: The  offset vectors that define the set \eqref{eq:27}. Right: Reconstruction of the four sources. The axes are as in Figure \ref{fig:3}.}
	\label{fig:4}
\end{figure}

In Figure \ref{fig:4} we show the reconstruction of $4$ sources. As in the previous example, the sources are located in the support of the CINT peak.
The image ${\cI}(\vz_o,\vy)$ shown in the top right plot has one spurious peak, due to the noise. However, this can be easily filtered 
out because the offset vector $\vec{\bf e} \in {\cal E}^{\rm{est}}$ corresponding to it does not have the property that $-\vec{\bf e} \in  {\cal E}^{\rm{est}}$. 
The set of remaining offset vectors is displayed in the left bottom plot and the reconstruction of the four sources is shown in the bottom right plot. The reconstruction is exact, up to the translation by $\vz_o$.

\section{Summary}
\label{sect:SUM}
We introduced a novel algorithm for array imaging in cluttered media modeled by a random sound speed. The algorithm is 
designed to work in the presence of strong scattering effects in clutter, where the sound waves recorded at the array 
are incoherent i.e., their statistical expectation is close to zero. Physically, this means that the range offset between the 
unknown sources and the array is larger than the scattering mean free path in clutter. The algorithm uses an imaging approach that 
is similar to the coherent interferometric (CINT) method. CINT is known to be robust to clutter scattering effects, as long as the 
waves are not in a diffusion regime i.e., for ranges less than a transport mean free path. The robustness 
comes at the  cost of image blur. This impedes imaging of sources at nearby locations $\vy_s$, for $s = 1, \ldots, N_s$. 
The algorithm introduced in this paper uses the observation that the blur affects only the estimation of the center locations 
$(\vy_s + \vy_{s'})/2$ of pairs $(s,s')$ of sources, whereas the offset vectors $\vy_s - \vy_{s'}$ can be estimated with the 
same  resolution as in the absence of clutter. Thus, it is possible to determine constellations of nearby sources, up to 
a translation within the support of a peak of the CINT image. 

We motivated the algorithm from first principles, starting with the wave equation in random media and assessed its performance with numerical simulations. To simplify the presentation, we considered a high frequency scattering regime defined by large, random wavefront distortions of the waves received at the array, although as explained in the paper, the results extend verbatim  to stronger scattering regimes.

\begin{appendices}
	
	\section{}
\label{sect:appA}
In this appendix we derive the expression \eqref{eq:kernel}.  We begin with the solution of the wave equation \eqref{eq:1}
evaluated at the receiver location $\vx_r = (\bx_r,0),$  
\begin{equation}
p(t,\vx_r) = \int_{-\infty}^\infty \frac{d \om}{2 \pi} e^{-i \om t} \hat f(\om) \sum_{s=1}^{N_s} \hat G(\om,\vx_r,\vy_s),
\label{eq:A1}
\end{equation}
where we used the assumption \eqref{eq:13} and $\hat G$ is the Green's function  modeled by \eqref{eq:15}. This model  holds under the   assumptions 
\begin{equation}
\ell \gg \sqrt{\la_oL} \gg \la_o, \quad \sigma \ll \frac{\sqrt{\la_o \ell}}{ L},
\label{eq:A8}
\end{equation}
as explained in \cite{RANDTT,CINT1,CINT_opt}. These ensure that the waves propagate 
along straight rays, as in geometrical optics, and that the random fluctuations of the amplitude of the Green's function 
are negligible. The first bound on $\ell$ is to have consistent assumptions on  $\sigma$,
\begin{equation}
\frac{\la_o}{\sqrt{\ell L}} \ll \frac{\la_o^{2/3} \ell^{1/6}}{L^{5/6}} \ll  \sigma \ll \frac{\sqrt{\la_o \ell}}{ L},
\label{eq:A9}
\end{equation}
chosen large enough to give large random travel time fluctuations, as explained in section \ref{sect:assume}.

It is shown in \cite[Proposition 3.1]{CINT_opt} that the expectation of the Green's function is 
\begin{equation}
\EE\left[\hat G(\om,\vx_r,\vy_s)\right] \approx \frac{e^{i \om \tau(\vx_r,\vy_s)}}{4 \pi |\vx_r-\vy_s|} e^{-\frac{|\vx_r-\vy_s|}{\cS}}.
\label{eq:A01}
\end{equation}
The decaying exponential is due to the random phase, which is approximately Gaussian, and the scattering mean free path is 
defined by 
\begin{equation}
\cS = \frac{8\la_o^2}{(2\pi)^{5/2} \sigma^2 \ell}.
\label{eq:A02}
\end{equation}
The lower bound on $\sigma$ in \eqref{eq:A8} implies that 
\begin{equation}
\cS \ll L = O(|\vx_r-\vy_s|),
\label{eq:A03}
\end{equation} 
so the wave recorded at the array is incoherent i.e., 
\begin{equation}
\EE[p(t,\vx_r)] \approx 0. 
\end{equation}

The kernel of the imaging function \eqref{eq:10} is  obtained from equations  \eqref{eq:6}, \eqref{eq:20} 
and \eqref{eq:A1},
\begin{align}
&{\cal K}(\vy,\vy',\vz,\vz') = \sum_{r,r'=1}^{N_r} \Psi \Big(\frac{|\bx_r-\bx_{r'}|}{X}\Big) \int_{-\infty}^{\infty} 
\frac{d \om}{2 \pi}  \times \nonumber \\
&\hspace{0.05in}\int_{-\infty}^{\infty}\frac{d \tom}{2 \pi}\, \hat \Phi \Big(\frac{\tom}{\Omega}\Big) \hat f \Big(\om + \frac{\tom}{2} \Big) 
\hat f^\star \Big(\om - \frac{\tom}{2} \Big) \hat G(\om,\vx_r,\vz) \times \nonumber \\
&\hspace{0.05in}\hat G^\star(\om,\vx_{r'},\vz') e^{-i \big(\om + \frac{\tom}{2}\big)\tau(\vx_r,\vy) + i \big(\om - \frac{\tom}{2}\big)\tau(\vx_{r'},\vy')},
\label{eq:A2}
\end{align}
where  we replaced $\vy_s$ by $\vz$ and $\vy_{s'}$ by $\vz'$, to avoid carrying over the source indexes. 
We use the definitions of $\hat f$, $\Psi$ and $\hat \Phi$ given  in section \ref{sect:assume}, and replace the sum over the receivers by the 
integral over the aperture, with Gaussian apodization,
\[
\sum_{r=1}^{N_r} ~\leadsto~ \frac{N_r}{a^2} \int_{\mathbb{R}^2} d \bx \, e^{-\frac{|\bx|^2}{2 (a/6)^2} }\, .
\]

\subsection{Wave decorrelation and the paraxial approximation}
Essentially the same calculation as in \cite[Section 4]{CINT1} shows that the kernel is statistically stable i.e., it is approximated by its expectation, 
when \[X = O(X_d) \ll a.\] The random travel time model accounts only for wavefront distortion and does not take into consideration 
delay spread due to scattering. Thus, the bandwidth does not play a big role in the statistical stability. However, in stronger scattering 
regimes the bandwidth is very important \cite{CINTTH} and statistical stability is achieved if  \[\Omega = O(\Omega_d) \ll B,\] as we assume here.  

The second moment formula is derived in \cite[Appendix B]{CINT_opt}
\begin{align}
&\EE\left[\hat G\Big(\om+\frac{\tom}{2},\vx,\vz\Big) \hat G^\star\Big(\om-\frac{\tom}{2},\vx',\vz'\Big)
\right] \nonumber \\
&\approx \frac{e^{i \om \widetilde{\tau}(\vx,\vx',\vz,\vz') + i \tom \overline{\tau}(\vx,\vx',\vz,\vz')}
}{(4 \pi)^2 |\vx-\vz| |\vx'-\vz'|} e^{-\frac{1}{\cS}| |\vx-\vz|-|\vx'-\vz'||}\times \nonumber \\
&\hspace{0.2in}e^{- \frac{1}{2 X_d^2}(|\bz'-\bz|^2 + (\bz'-\bz) \cdot (\bx'-\bx) + |\bx'-\bx|^2)- \frac{\tom^2}{2 \Omega_d^2}},
\label{eq:A3}
\end{align}
with $\overline\tau$ and $\widetilde{\tau}$ defined in \eqref{eq:defOTTau1}--\eqref{eq:defOTTau2}.
It decays  with  the frequency and cross-range  offsets, due to the decorrelation of the waves in the 
random medium. The decorrelation frequency  and length are 
\begin{equation}
\Omega_d = \frac{2 \om_o}{(2 \pi)^{5/4}} \left(\frac{\la_o}{\sigma \sqrt{\ell L}}\right) \ll \om_o, 
\label{eq:A4}
\quad  X_d = \sqrt{3} \ell \frac{\Omega_d}{\om_o} \ll \ell,
\end{equation}
where the inequalities  are implied by \eqref{eq:A8}.

As stated in section \ref{sect:assume}, we consider a paraxial wave propagation regime, where 
\begin{align}
\om_o \tau(\vx,\vz) = k_o |\vx-\vz| \approx k_o \Big( z_3 + \frac{|\bx-\bz|^2}{2 L} \Big), 
\label{eq:A6}
\end{align} 
with negligible residual 
\begin{equation}
O\Big(\frac{a^4}{\la_o L^3}\Big) + 
O\Big(\frac{a^2 D_3}{\la_o L^2}\Big) \ll 1.
\label{eq:A7}
\end{equation}
Here we used the scaling relation \eqref{eq:17p}. We also approximate the amplitude of the Green's functions by 
\begin{equation}
\frac{1}{4 \pi |\vx-\vz|} \approx \frac{1}{4 \pi L}.
\end{equation}

Since the expression \eqref{eq:A3} is large when the cross-range offsets are $O(X_d)$, we estimate from \eqref{eq:A02},  
\eqref{eq:A8} and \eqref{eq:A7}  that 
\begin{equation}
\frac{||\vx-\vz|-|\vx'-\vz'||}{\cS}  \ll 1,
\label{eq:A10}
\end{equation}
so the decaying exponential in the second line of \eqref{eq:A3} is approximately equal to $1$. 

\subsection{Calculation of the imaging kernel}
Because of the decorrelation of the waves over the scale $X_d$, it suffices to consider imaging points at cross-range offsets 
$|\by-\by'| \le O(X_d)$. 
Substituting the results  in \eqref{eq:A2}, and using $\Omega_e$ and $X_e$ defined in \eqref{eq:23}--\eqref{eq:24}, 
which are similar to $\Omega_d$ and $X_d$, we obtain that 
\begin{align}
{\cal K}&(\vy,\vy',\vz,\vz') \approx \frac{N_r^2}{2 \pi (4 \pi L)^2 B a^4} e^{-\frac{|\tbz|^2}{2 X_d^2}} \times \nonumber \\
&\int_{\mathbb{R}^2} d \bx \, 
e^{-\frac{|\bx|^2}{(a/6)^2}} \int_{\mathbb{R}^2} d \tbx \, e^{-\frac{|\tbx|^2}{2 X_e^2} - \frac{\tbx \cdot \tbz}{2 X_d^2}}
\times \nonumber \\
&\int_{-\infty}^\infty d \om \, e^{-\frac{(\om-\om_o)^2}{2 B^2} + i 
\frac{\om}{c_o} \big[ \tz_3-\ty_3  - \frac{\bx \cdot (\tbz-\tby)}{L} - \frac{\tbx \cdot (\cbz-\cby)}{L} + 
\frac{\cbz \cdot \tbz - \cby\cdot \tby}{L} \big]} \times \nonumber \\
&\int_{-\infty}^\infty d \tom \, e^{-\frac{\tom^2}{2 \Omega_e^2}+ i \frac{\tom}{c_o} 
\big[ \cz_3 - \cy_3 - \frac{\bx \cdot (\cbz-\cby)}{L} + \frac{|\cbz|^2-|\cby|^2}{2 L}\big] },
\label{eq:A11}
\end{align}
with center  vectors $(\cby,\cy_3)$ and difference vectors $(\tby,\ty_3)$ defined in \eqref{eq:21}--\eqref{eq:22}. The vectors 
$(\cbz,\cz_3)$ and  $(\tbz,\tz_3)$ are defined the same way, by replacing $\vy$ and $\vy'$ with $\vz$ and $\vz'$ in \eqref{eq:21}--\eqref{eq:22}.
Note that in \eqref{eq:A11} we neglect the phase terms 
\[
\frac{\tom}{c_o} \left(\frac{|\tbz|^2-|\tby|^2}{4L} \right) = O\left( \frac{\Omega_e X_d^2}{c_o L} \right) \ll 1 
\]
and 
\[
\frac{\tom}{c_o} \frac{\tbx \cdot (\tbz-\tby )}{L} = O\left( \frac{\Omega_e X_d^2}{c_o L} \right) \ll 1,
\]
with the inequalities implied by \eqref{eq:25}, \eqref{eq:A8} and \eqref{eq:A4}. 

Carrying out the Gaussian integrals  in \eqref{eq:A11}, we obtain 
\begin{align}
{\cal K}(\vy,\vy',&\vz,\vz') \approx H \exp \Big\{-\frac{|\tbz|^2}{2 \gamma X_d^2} - \frac{\beta^2}{2 \big(1 + \frac{|\czet|^2}{2 \theta^2}\big)}
  \nonumber\\
&-\frac{\upsilon^2 + \eta^2+ 2i \eta \upsilon^2 B/\om_o}{2\big(1 + \frac{\upsilon^2 B^2}{\om_o^2}\big)}  + i \frac{\om_o \eta}{B}\Big\}, 
\label{eq:A20}
\end{align}
where the amplitude factor is 
\begin{align*}
H &= \frac{N_r^2 X_e^2 \Omega_e}{288 a^2 L^2 \big(1 + \frac{|\czet|^2}{2 \theta^2}\big)^{1/2} 
\big(1 + \frac{\upsilon^2 B^2}{\om_o^2} \big)^{1/2}}, 
\end{align*}
and we used the positive constants $\theta$ and $\gamma$ defined by 
\[\theta = \frac{\om_o X_e}{\Omega_e (a/6)}, \quad \frac{1}{\gamma} = 1-\frac{X_e^2}{4 X_d^2}  > \frac{3}{4}. \]
The inequality on $\gamma$ is implied by definition \eqref{eq:24}. We also  introduced the  notation 
\begin{align*} 
\czet &= \frac{\cbz-\cby}{L/(k_o X_e)}, \quad \tzet = \frac{\tbz-\tby}{6 \sqrt{2} L/(k_o a)}, \\
\beta &= \frac{\cz_3 - \cy_3 + \frac{|\cbz|^2 - |\cby|^2}{2 L}}{c_o/\Omega_e} \approx \frac{|(\cbz,\cz_3)|-|(\cby,\cy_3)|}{c_o/\Omega_e}, \\
\upsilon^2 &= |\czet|^2 + \frac{|\tzet|^2}{1+ \frac{|\czet|^2}{2 \theta^2}} + \frac{|\czet|^2 |\tzet|^2 - |\czet \cdot \tzet|^2}{
	 { {2 \theta^2} \left( 1+ \frac{|\czet|^2}{2 \theta^2} \right)  
}
 }, 
\end{align*}
and 
\begin{align*}
\eta = \frac{B}{c_o} \left[ \tz_3 - \ty_3 + \frac{\cbz \cdot \tbz - \cby \cdot \tby}{L} +  \frac{X_e^2}{2 X_d^2} \frac{\tbz \cdot (\cbz-\cby)}{L} \right. \\ 
\left.- 
\frac{\beta \czet \cdot \tzet}{ \sqrt{2} k_o \theta \big( 1 + \frac{|\czet|^2}{2 \theta^2}\big)}\right].
\end{align*}

The expression \eqref{eq:A20} simplifies, because it is large only when $|\upsilon|, |\beta|, |\eta| = O(1)$. Since $B \ll \om_o$, 
we can write
\[
1 + \frac{\upsilon^2 B^2}{\om_o^2} \approx 1,
\]
and neglect the $\eta \upsilon^2 B/\om_o$ phase,  with absolute value much less than $1$. We also have 
that $|\czet|$ and $|\tzet|$ are $O(1)$, because $\upsilon = O(1)$, and therefore 
\begin{align*}
\frac{B}{c_o} &\left|\frac{\cbz \cdot \tbz - \cby \cdot \tby}{L} \right| = 
\frac{B}{c_o} \left| \frac{\cbz \cdot (\tbz-\tby) + \tby \cdot (\cbz-\cby)}{L} \right| \\
&= \frac{B}{\om_o} \left[O \left( \frac{D k_o |\tbz-\tby|}{L} \right) + O\left(\frac{X k_o |\cbz-\cby|}{L} \right) \right]\\
&= \frac{B}{\om_o} \left[ O\left(\frac{D}{a}\right) + O\left(\frac{X}{X_e}\right) \right] \ll 1.
\end{align*}
Similarly, we obtain that 
\begin{align*}
\frac{B}{c_o}  \left|\frac{X_e^2}{2 X_d^2} \frac{\tbz \cdot (\cbz-\cby)}{L}\right| \ll 1,
\end{align*}
and 
\begin{align*}
\frac{B}{c_o} \left| \frac{\beta \czet \cdot \tzet}{k_o \theta \big( 1 + \frac{|\czet|^2}{2 \theta^2}\big)} \right| = O \left(\frac{B}{\om_o \theta}\right).
\end{align*}
The definitions of $\theta$, $X_e$ and $\Omega_e$  give 
\[
\theta = O\left(\frac{\ell}{a} \right) \ll 1.
\]
We assume that this number is not too small, so that $B/(\om_o \theta) \ll 1$,  and we can use the results above to approximate 
\[
\eta \approx \frac{B}{c_o} (\tz_3-\ty_3).
\]

We can now rewrite the kernel \eqref{eq:A20} as 
\begin{align} 
{\cal K}(\vy,\vy',\vz,\vz') &\approx H \exp \Big\{ i \frac{\om_o \eta}{B} - 
\frac{|\tbz|^2}{2 \gamma X_d^2} - \frac{(\tz_3-\ty_3)^2}{2 (c_o/B)^2}
  \nonumber\\
&-\frac{\big| |(\cbz,\cz_3)| - |(\cby,\cy_3| \big|^2}{2 (c_o/\Omega_e)^2\big(1 + \frac{|\czet|^2}{2 \theta^2}\big)}-\frac{\upsilon^2 }{2}\Big\}. 
\label{eq:A21}
\end{align}
Taking absolute value and using that 
\begin{align*}
\upsilon^2 \ge  |\czet|^2 + \frac{|\tzet|^2}{1 + \frac{|\czet|^2}{2 \theta^2}},
\end{align*}
we obtain the result \eqref{eq:kernel} with
 $ \gamma_1 =   6\sqrt{2 \left(1 + \frac{|\czet|^2}{2 \theta^2} \right)} \ge O(1).$

\section{}
\label{sect:appB}
{Suppose that there exists a constellation $\mathcal{Y}_{0}$
of sources such that $\mathcal{E}^{est}=\mathcal{E}\left(\mathcal{Y}_{0}\right)$.
We show here that 
	\begin{equation}
	\mathcal{Y}^{est}=\textbf{Algorithm}(\mathcal{E}^{est},\left\{ \vz_{0}\right\} )\label{MainCall}
	\end{equation}
returns a set $\mathcal{Y}^{est}$ such that 
\begin{equation}
\mathcal{E}\left(\mathcal{Y}^{est}\right)=\mathcal{E}^{est}.
\label{eq:success}
\end{equation}

Once we show that $\mathcal{Y}^{est} \ne \emptyset$, it is straightforward to see from the 
definition of \textbf{Algorithm} that \eqref{eq:success} must holds. It remains to show that $\mathcal{Y}^{est}$ is not an empty set.

For a proof by contradiction, suppose that  \eqref{MainCall} returns $\mathcal{Y}^{est} = \emptyset$. This means explicitly  that each call of $\textbf{Algorithm}$ results in executing line 14.
Since the set $\mathcal{E}({\cal Y})$ of offsets is translation invariant, let us replace 
the set ${\cal Y} = \{\vy_1, \ldots, \vy_{n_s}\}$ of source locations  by the set ${\cal Y}_o = \{\vz_{0},\vz_{1},\vz_{2},\dots,,\vz_{n_s-1}\} $ of translated source locations, with translation defined by $\vz_o$.  
Define the vectors
\begin{equation}
\label{eq:veces}
\vec{{\bf e}}_{j_{k}}=\vz_{k}-\vz_{0},\quad k = 1, \ldots, m, ~ ~ m = n_s-1,
\end{equation}
which belong  to $\mathcal{E}^{est}$. Note that $\mathcal{E}^{est}$ contains other offset vectors, as well.
Without loss of generality, we can assume that the vectors \eqref{eq:veces} are enumerated in order in $\mathcal{E}^{est}$,
meaning that $\vec{{\bf e}}_{j_{k}}$ comes
before $\vec{{\bf e}}_{j_{k+1}}$ for $k\in\left\{ 1,\dots,m-1\right\} $.
\vspace{0.1in}
\begin{description}
\itemsep 0.1in
	\item [$1^{\rm{st}}$-Recursion:] In the first call of \textbf{Algorithm}, the
	arguments are $E=\mathcal{E}^{est}$ and $Y=\left\{ \vz_{0}\right\} $.
	Since we assume \eqref{MainCall} returned $\emptyset$,
	the line 14 was executed. Thus in the while loop, each element of $E$ is selected
	at line 5. In particular, at  line
	5, the vector $\vec{{\bf e}}=\vec{{\bf e}}_{j_{1}}$
	is removed from the set $E$, and at line 6 we have 
	\[
	\vy=\vz_{0}+\vec{{\bf e}}=\vz_{0}+\vec{{\bf e}}_{j_{1}}=\vz_{1}.
	\]
	But then 
	\[
	\left\{ \pm\left(\left\{ \vz_{0}\right\} -\vy\right)\right\} =\left\{ \pm\vec{{\bf e}}_{j_{1}}\right\} \subseteq\mathcal{E}^{est}.
	\]
	Therefore line 8 has to be executed. This will take us to the next recursion.
	\item [{$2^{\rm{nd}}-$Recursion:}] At this recursion level, $E$ contains the offset vectors 
	$\vec{{\bf e}}_{j_{2}},\dots,\vec{{\bf e}}_{j_{m}}$,
	and $Y=\left\{ \vz_{0},\vz_{1}\right\} $. We note that $Y\subset\mathcal{Y}_{0}$ at every recursion level.
	
	Once again, by assumption, line 14 was executed. Therefore,
	at some iteration of the while loop, at line 5, the vector $\vec{{\bf e}}=\vec{{\bf e}}_{j_{2}}$
	is removed from the set $E$, and in line 6, 
	\[
	\vy=\vz_{0}+\vec{{\bf e}}=\vz_{0}+\vec{{\bf e}}_{j_{2}}=\vz_{2}.
	\]
	We have 
	\[
	\left\{ \pm\left(Y-\vy\right)\right\}  \subseteq\left\{ \pm\left(\mathcal{Y}_{0}-\vz_{2}\right)\right\} \subseteq\mathcal{E}^{est}
	\]
	since $Y\subset\mathcal{Y}_{0}$, so line 8 must  be executed. This takes us to the next recursion, 
	with $E$ containing the offset vectors 
	$\vec{{\bf e}}_{j_{3}},\dots,\vec{{\bf e}}_{j_{m}}$,
	and $Y\cup\left\{ \vy\right\} =\left\{ \vz_{0},\vz_{1},\vz_{2}\right\} $.
	\item [{$m^{\rm{th}}-$Recursion:}] At this recursion level, $E$ contains the
	element $\vec{{\bf e}}_{j_{m}}$ and $Y=\left\{ \vz_{0},\vz_{1},\dots,\vz_{m-1}\right\} $.
	Similar to above, we know that at some iteration of the while
	loop, at line 5, the vector $\vec{{\bf e}}=\vec{{\bf e}}_{j_{m}}$
	is removed from the $E$ , and in line 6, 
	\[
	\vy=\vz_{0}+\vec{{\bf e}}=\vz_{0}+\vec{{\bf e}}_{j_{m}}=\vz_{m}.
	\]
	Moreover 
	\[
	\left\{ \pm\left(Y-\vy\right)\right\}\subseteq\left\{ \pm\left(\mathcal{Y}_{0}-\vz_{m}\right)\right\} \subseteq\mathcal{E}^{est}
	\]
	since $Y\subset\mathcal{Y}_{0}$. Therefore line 8 has to be executed,
	and this  take us to next recursion with all the vectors \eqref{eq:veces} removed from $E$ and $Y\cup\left\{ \vy\right\} =\left\{ \vz_{0},\vz_{1},\dots,\vz_{m-1},\vz_{m}\right\} $. 
	\item [{$\left(m+1\right)^{\rm{th}}-$Recursion:}] At this recursion level since 
	\[Y=\left\{ \vz_{0},\vz_{1},\dots,\vz_{m}\right\}  = {\cal Y}_o\]
	and since 
	$
	\mathcal{E}^{est}=\mathcal{E}\left({\cal Y}_o\right)$, line 
	2 is executed and the recursion returns $\mathcal{Y}^{est}=Y$, which is a non-empty set.
\end{description}
This contradicts  our assumption
$\mathcal{Y}^{est}=\emptyset$. Thus, we conclude \eqref{MainCall} returns
a non-empty $\mathcal{Y}^{est}.$

}

\end{appendices}

\section*{Acknowledgment}
This work is supported in part by the NSF grant DMS1510429 and by the NSF Grant  DMS-1439786, while Liliana Borcea was in residence at the Institute for Computational and Experimental Research in Mathematics in Providence, RI, during the Fall 2017 semester.

%



\end{document}

%% file: mathmacros.tex
\renewcommand{\hat}{\widehat}

\newcommand{\om}{\omega}

\newcommand{\la}{\lambda}

\newcommand{\cA}{\mathcal{A}}
\newcommand{\cS}{\mathcal{S}}

\newcommand{\EE}{\mathbb{E}}
\newcommand{\cI}{\mathcal{I}}

\newcommand{\vy}{{\vec{\bm{y}}}}
\newcommand{\vx}{\vec{\bm{x}}}
\newcommand{\vz}{\vec{\bm{z}}}
\newcommand{\by}{\bm{y}}
\newcommand{\bz}{\bm{z}}
\newcommand{\bn}{{\bm{n}}}
\newcommand{\bx}{\bm{x}}

\newcommand{\tom}{\widetilde \omega}
\newcommand{\tbx}{\widetilde {\bx}}
\newcommand{\tby}{\widetilde {\by}}
\newcommand{\ty}{\widetilde {y}}
\newcommand{\cy}{\overline {y}}
\newcommand{\tbz}{\widetilde{\bz}}
\newcommand{\cbz}{\overline{\bz}}
\newcommand{\cby}{\overline{\by}}
\newcommand{\tz}{\widetilde{z}}
\newcommand{\cz}{\overline{z}}
\newcommand{\tit}{\widetilde{t}}
\newcommand{\czet}{\overline{\boldsymbol{\zeta}}}
\newcommand{\tzet}{\widetilde{\boldsymbol{\zeta}}}

\makeatletter
\newsavebox\myboxA
\newsavebox\myboxB
\newlength\mylenA

\newcommand*\xoverline[2][1.3]{%
    \sbox{\myboxA}{$\m@th#2$}%
    \setbox\myboxB\null
    \ht\myboxB=\ht\myboxA%
    \dp\myboxB=\dp\myboxA%
    \wd\myboxB=#1\wd\myboxA
    \sbox\myboxB{$\m@th\overline{\copy\myboxB}$}
    \setlength\mylenA{\the\wd\myboxA}
    \addtolength\mylenA{-\the\wd\myboxB}%
    \ifdim\wd\myboxB<\wd\myboxA%
       \rlap{\hskip 0.3\mylenA\usebox\myboxB}{\usebox\myboxA}%
    \else
        \hskip -0.3\mylenA\rlap{\usebox\myboxA}{\hskip 0.5\mylenA\usebox\myboxB}%
    \fi}
\makeatother